\documentclass[11pt]{amsart}

\usepackage[english]{babel}
\usepackage[T1]{fontenc}
\usepackage[utf8]{inputenc}

\usepackage{mathpazo}
\usepackage{graphicx}
\usepackage{microtype}
\usepackage{enumerate}
\usepackage{csquotes}
\usepackage[dvipsnames]{xcolor}
\usepackage[bookmarks, colorlinks, citecolor=MidnightBlue!75!black, urlcolor=MidnightBlue!75!black, linkcolor=MidnightBlue!75!black]{hyperref}

\calclayout

\AtBeginDocument{%
   \def\MR#1{}
}

\usepackage{amsmath, amssymb, amsthm, amscd}
\usepackage{commath, mathtools, mathrsfs, thmtools}
\usepackage{tikz, tikz-cd}
\usetikzlibrary{matrix, arrows, decorations.pathmorphing}

\theoremstyle{plain}
\newtheorem{theorem}{Theorem}[section]
\newtheorem{lemma}[theorem]{Lemma}
\newtheorem{proposition}[theorem]{Proposition}
\newtheorem{corollary}[theorem]{Corollary}

\newtheorem*{theorem*}{Theorem}
\newtheorem*{lemma*}{Lemma}
\newtheorem*{proposition*}{Proposition}
\newtheorem*{corollary*}{Corollary}

\newtheorem*{TA}{Theorem A}
\newtheorem*{TB}{Theorem B}
\newtheorem*{TC}{Theorem C}
\newtheorem*{TD}{Theorem D}

\theoremstyle{definition}
\newtheorem{remark}[theorem]{Remark}
\newtheorem{definition}[theorem]{Definition}
\newtheorem{example}[theorem]{Example}

\newtheorem*{conjecture*}{Conjecture}
\newtheorem*{remark*}{Remark}
\newtheorem*{definition*}{Definition}
\newtheorem*{observation*}{Observation}

\newcommand{\on}{\operatorname}

\newcommand{\A}{\mathbb{A}}

\newcommand{\aut}{\on{Aut}}

\newcommand{\C}{\mathbb{C}}
\newcommand{\ced}{\on{ced}}

\newcommand{\ed}{\on{ed}}
\newcommand{\et}{\textrm{ét}}
\newcommand{\fced}{\on{fced}}
\newcommand{\fed}{\on{fed}}

\newcommand{\fld}{\on{Fields}}

\newcommand{\gal}{\on{Gal}}

\newcommand{\G}{\mathbb{G}}

\newcommand{\hz}{\hat{\mathbb{Z}}}
\newcommand{\id}{\on{id}}

\newcommand{\mb}{\mathbb}
\newcommand{\mc}{\mathcal}

\newcommand{\N}{\mathbb{N}}

\newcommand{\Q}{\mathbb{Q}}

\newcommand{\s}{\subseteq}

\newcommand{\spec}{\on{Spec}}

\newcommand{\trdeg}{\on{trdeg}}
\newcommand{\uaut}{\underline{\on{Aut}}}
\newcommand{\uisom}{\underline{\on{Isom}}}
\newcommand{\upi}{\underline{\pi}}

\newcommand{\xar}{\xrightarrow}
\newcommand{\Z}{\mathbb{Z}}

\renewcommand{\epsilon}{\varepsilon}
\renewcommand{\H}{\on{H}}
\renewcommand{\hom}{\on{Hom}}
\renewcommand{\injlim}{\varinjlim}

\renewcommand{\O}{\mc{O}}

\renewcommand{\phi}{\varphi}
\renewcommand{\projlim}{\varprojlim}

\renewcommand{\tilde}{\widetilde}
\renewcommand{\hat}{\widehat}

\usepackage[giveninits=true,style=numeric,backend=biber,maxnames=99,maxalphanames=5,giveninits=true]{biblatex}
\bibliography{fed}

\author{Giulio Bresciani}
\address{Freie Universit\"at Berlin, Arnimallee 3, 14195, Berlin, Germany}
\email{gbresciani@zedat.fu-berlin.de}
\title{Essential dimension and pro-finite group schemes}
\date{}

\subjclass[2010]{14L15, 14F35, 11G35}
\thanks{The author is supported by the DFG Priority Program "Homotopy Theory and Algebraic Geometry" SPP 1786}

\begin{document}

\begin{abstract}
	A. Vistoli observed that, if Grothendieck's section conjecture is true and $X$ is a smooth hyperbolic curve over a field finitely generated over $\Q$, then $\upi_{1}(X)$ should somehow have essential dimension $1$. We prove that an infinite, pro-finite étale group scheme always has infinite essential dimension. We introduce a variant of essential dimension, the fce dimension $\fced G$ of a pro-finite group scheme $G$, which naturally coincides with $\ed G$ if $G$ is finite but has a better behaviour in the pro-finite case. Grothendieck's section conjecture implies $\fced\upi_{1}(X)=\dim X=1$ for $X$ as above. We prove that, if $A$ is an abelian variety over a field finitely generated over $\Q$, then $\fced\upi_{1}(A)=\fced TA=\dim A$.
\end{abstract}

\maketitle
\tableofcontents

\section{Introduction}

Essential dimension has been introduced by Buhler and Reichstein in \cite{br97} as a measure of the complexity of torsors under a group scheme $G$: it is the minimal number of parameters necessary in order to define a generic $G$-torsor. Up to now, essential dimension has only been studied for group schemes of finite type, for which a number of tools has been developed: most notably, the existence of the so-called versal torsors, which among other things ensures that every affine group scheme of finite type has finite essential dimension.

\subsection{Essential dimension of pro-finite group schemes}
The situation is completely different for pro-finite group schemes, but in some sense much simpler.

\begin{TA}[\ref{crit}]
	Let $G$ be a pro-finite, infinite étale group scheme over a field $k$. Then there exists a $G$-torsor of infinite essential dimension.
	
	In particular, a pro-finite étale group scheme is finite if and only if it has finite essential dimension.
\end{TA}

Observe that Theorem A does not only say that $\ed G=\infty$: a priori, it might happen that $G$-torsors have finite, arbitrarily large essential dimension. We are able to construct a single torsor of infinite essential dimension.

Because of Theorem A, one can either stop studying essential dimension of pro-finite group schemes or try to modify its definition in order to get more interesting results.

\subsection{Essential dimension and anabelian geometry}

We take the second path, mainly because of the following observation due to A. Vistoli.

\begin{observation*}
Let $X$ be a smooth, projective hyperbolic curve of genus at least $2$ with a rational point $x$ over a field $k$ finitely generated over $\Q$. If $\upi_{1}(X,x)$ is the étale fundamental group scheme of $X$, Grothendieck's section conjecture predicts that the natural map $X(k')\to\H^{1}(k',\upi_{1}(X,x))$ is bijective for every finitely generated extension. Thus, $\upi_{1}(X,x)$ should somehow have essential dimension equal to $\dim X=1$.
\end{observation*}

The reason why Vistoli's observation fails if we interpret it in a na\"ive sense is that essential dimension depends on all extensions of the base fields, while the section conjecture is only expected to be true for fields finitely generated over $\Q$. We call \emph{finite type essential dimension} $\fed_{k} G$ of a pro-finite group scheme $G$ the supremum of essential dimensions of $G$-torsors defined over fields finitely generated over $k$, $\fed G=\ed G$ is $G$ is finite. Grothendieck's section conjecture implies $\fed\upi_{1}(X)=\dim X=1$. 

There is also a version of Grothendieck's conjecture for non-projective curves, but if $X$ is not projective and $X\neq\A^{1}$ we show that $\fed\upi_{1}(X)=\infty$. This happens because of "pathological" phenomena appearing when passing to the projective limit without considering pro-finite topologies, similarly to how na\"ive Galois theory for infinite extensions is just false. Taking into account these problems leads to a further refinement of essential dimension, the \emph{fce dimension} $\fced G$ of a pro-finite group scheme $G$. Grothendieck's section conjecture implies  $\fced\upi_{1}(X)=\dim X=1$, both in the affine and projective case. 

It might be possible to compute $\fced\upi_{1}(X)=\dim X=1$ without knowing the section conjecture. As we will see later in the introduction, we obtain a result that goes in this direction: if $A$ is an abelian variety over a field finitely generated over $\Q$, then $\fced\upi_{1}(A)=\dim A$.

If one could prove $\fced\upi_{1}(X)=\dim X=1$, we do not expect it to directly help prove the section conjecture, but it would give evidence: it is a non-trivial test.

In this paper, we start to compute these variants of essential dimension for several instances of abelian pro-finite group schemes, since the subject is still unexplored.

\subsection{Finite type essential dimension}

We obtain a complete description of finite type essential dimension for abelian, torsion free pro-$l$ group schemes.

If $A$ is an abelian group scheme and $l$ is any prime number, denote by $T_{l}A=\projlim_{n}A[l^{n}]$ its $l$-adic Tate module, it is a (possibly ramified) pro-$l$ group scheme over $k$. Moreover, write $TA=\projlim_{n}A[n]=\prod_{l}T_{l}A$ for the global Tate module.

\begin{definition*}
	Let $l$ be a prime and $G$ a pro-$l$ group scheme over a field $k$ with $\on{char}k\neq l$. We say that $G$ \emph{has a geometric subgroup} if there exists a finite extension $k'/k$ and a semi-abelian variety $A$ over $k'$ with a non-trivial homomorphism $T_{l}A\to G_{k'}$.
\end{definition*}

\begin{TB}[\ref{fed}]
	Let $l$ be a prime and $G$ an abelian, torsion free, pro-$l$ group scheme over a field $k$ with $\on{char}k\neq l$. If $G$ has a geometric subgroup, then $\fed G=\infty$, otherwise $\fed G=0$.
\end{TB}

If $k$ is finitely generated over $\Q$, Deligne's theory of weights gives an obstruction to the existence of a geometric subgroup.

\begin{corollary*}[\ref{fedwp}]
	Let $k$ be a field finitely generated over $\Q$, and let $l$ be a prime number. Let $G$ be an abelian, torsion free pro-finite $l$-group scheme over $k$. Assume that the Galois representation $G\otimes_{\Z_{l}}\Q_{l}$ is pure of weight $w\neq-2,-1$. Then $\fed_{k} G=0$.
\end{corollary*}

Theorem B shows that finite type essential dimension depends heavily on the arithmetic of the base field: for instance, $\fed_{k}\Z_{l}=0$ if $k$ is a number field but $\fed_{k}\Z_{l}=\infty$ if $k=\overline{\Q}$, because $\Z_{l}\simeq\Z_{l}(1)$ over $\overline{\Q}$.

A theorem of Florence \cite[Theorem 4.1]{flo08} implies that $\lim_n\ed_{k}\Z/l^{n}=\infty$ over fields finitely generated over $\Q$, but still we have $\fed_{k}\Z_{l}=0$. On the other hand, we have that $\ed_{k}\mu_{l^{n}}=1$ for every $n$ but $\fed_k\Z_{l}(1)=\infty$.

Let us also mention that for a non-trivial abelian variety $A$ over a number field $k$, Brosnan and Sreekantan have shown in \cite{bs08} that $\ed_k A=\infty$, and they do so by showing that, for some prime $l$, $\ed_k A[l^{n}]$ diverges. Note, however, that there is little connection between their result and Theorem B for $T_{l}A$: our result is valid for every prime different from the characteristic and over an arbitrary base field, where as the theorem of Brosnan and Sreekantan holds only for a fixed prime and over number fields, the proofs are completely different. Over $\C$, Brosnan has shown that $\ed_{\C}A=2\dim A$, see \cite{bro07}.

\subsection{Continuous and fce dimensions}

The proof that $\fed G=\infty$ when $G$ has a geometric subgroup relies on a pathological construction: it is analogous to the fact that $\Z_{p}$ has infinite rank as a $\Z$-module, even if its "topological rank" is finite. 

We can solve this problem by introducing a "topological" variant of essential dimension for pro-finite groups (or gerbes) which naturally coincides with the classical one for finite groups, see \ref{ceddef} for the definition. We call it \emph{continuous essential dimension} $\ced G$ of a group $G$. One can also merge the finite and the continuous variants, obtaining what we call the \emph{fce dimension} $\fced G$ of $G$. 

We compute the continuous and fce dimension for Tate modules of tori (using flat cohomology in order to handle inseparability issues). Our proof works for torsors under tori, where the Tate module is replaced by Nori's fundamental gerbe $\Pi_{\bullet/k}$, see \cite{bv15} for details.

\begin{TC}[\ref{fcedzp1}]
	Let $D^{1}/k$ be a $D$-torsor over a field $k$, where $D/k$ is an algebraic torus and $l$ any prime number (possibly equal to $\on{char}k$). Then 
	\[\ced_{k}\Pi_{D^1/k,l}=\fced_k\Pi_{D^1/k,l}=\ced_{k}\Pi_{D^1/k}=\fced_k\Pi_{D^1/k}=\dim D^{1}.\]
	In particular, for $D^{1}=D$ we have
	\[\ced_{k}T_{l}D=\fced_{k}T_{l}D=\ced_{k}TD=\fced_{k}TD=\dim D.\]
\end{TC}

The proof of Theorem C is an easy application of Kummer theory plus Hilbert's theorem 90, but we think the result is interesting: the dimension of a torus $D$ coincides with the (finite) continuous essential dimension of its fundamental group. Anabelian geometry studies how geometric properties of a variety reflect in its fundamental group, and this is what is happening in Theorem C, even if in a shallow sense. Note that it is possible to find a torus $D$ with $\ed D[n]>\dim D$ for every $n$ divisible enough, see the \autoref{scavia} for which I thank F. Scavia: Theorem C is not just a consequence of the fact that $\on{rank}T_{p}D=\dim D$.

For abelian varieties, we cannot expect a result analogous to Theorem C to hold over any field: for instance, if $A$ is an abelian variety over $\C$, then $T_{l}A=\Z_{l}^{2\dim A}=\Z_{l}(1)^{2\dim A}$ and thus $\fced T_{l}A=2\dim A$ thanks to Theorem C. We prove that $\fced TA=\dim A$ holds for abelian varieties if the base field is finitely generated over $\Q$. The proof relies on Faltings' theorem (Tate conjecture for abelian varieties).

\begin{TD}[\ref{dscab}]
	Let $A^1$ be an $A$-torsor for an abelian variety $A$ over a field $k$ finitely generated over $\Q$, and $l$ a prime number. Then 
	\[\fced_k\Pi_{A^1/k}=\fced_k\Pi_{A^1/k,l}=\dim A^1.\]
	In particular, for $A^1=A$ we have
	\[\fced_kTA=\fced_kT_lA=\dim A.\]
\end{TD}

We really need the fce dimension in this case: by Theorem B we have $\fed_{k}TA=\infty>\dim A$, while by Theorem C
\[\ced_{k} TA\ge\ced_{\bar{k}}TA_{\bar{k}}=\ced_{\bar{k}}\hz^{2\dim A}=2\dim A>\dim A,\]
in contrast to what happens for tori.

\subsection{Extension of torsors}

An important problem when studying finite type essential dimension (either continuous or not) is deciding whether a torsor on the generic point of a variety extends to the whole variety. We study this problem in \autoref{extapp}, in particular we find a sufficient condition of arithmetic flavour, see \autoref{r-proper}. We stress that this study is a crucial and not merely technical ingredient in the proofs of both Theorem B and Theorem D.

\subsection*{Acknowledgements}
I would like to thank Angelo Vistoli for many useful discussions and for the proof of \autoref{proftor}, Federico Scavia for \autoref{scavia} and Zinovy Reichstein for uncountably many useful comments. Finally, I would also like to thank an anonymous referee for suggesting the use of Deligne weights.

\section{Notations and preliminaries}

A variety is a geometrically integral scheme of finite type over a field. A morphism $X\to Y$ of fibered categories (e.g. schemes or gerbes) over a field $k$ is \emph{constant} if there exists a factorization $X\to\spec k\to Y$.

If $X$ is an inflexible fibered category over $k$ (for instance, a geometrically connected, geometrically reduced scheme) we denote by $\Pi_{X/k}$ its Nori fundamental gerbe and by $\Pi_{X/k}^{\et}$ its étale fundamental gerbe, see \cite{bv15}. If $\on{char}k=0$, then $\Pi_{X/k}=\Pi^{\et}_{X/k}$. 

By an abuse of notation, we do not distinguish between (pro-)finite groups and the (pro-)discrete group schemes associated to them.

We will use extensively basic inequalities about essential dimension without explicit reference, see \cite{bf03} for a compendium. These include, but are not limited to,
\begin{itemize}
	\item if $T\to\spec K$ is a $G$-torsor and $K'/K/k$, $K'/k'/k$ are extensions, then $\ed_{k'}T_{K'}\le\ed_{k}T$,
	\item if $G$ is a group scheme over a field $k$ and $k'/k$ is an extension, then $\ed_{k'}G_{k'}\le\ed_{k}G$,
	\item if $T$ is a $G$-torsor and $G\to G'$ is an homomorphism, then $\ed_{k}T\times^{G}G'\le\ed_{k}T$.
\end{itemize}

If there is no risk of confusion, we will often drop the subscript $\bullet_{k}$, for instance we will write $\ed G$ instead of $\ed_{k}G$.

\subsection{$\H^1$ for pro-finite group schemes} If $G=\projlim_i G_i$ is a pro-algebraic group scheme over $k$ (every affine group scheme is pro-algebraic) and $L/k$ is an extension, we have to clarify what we mean by $\H^1(L,G)$: there are at least two possibilities.
\begin{itemize}
	\item $\projlim_i\H^1(L,G_i)$.
	\item The set $\on{Tors}(L,G)$ of $G$-torsors over $L$, where by $G$-torsor we mean a scheme $T$ over $L$ with an action of $G$ such that $G_L\times_L T\to T\times_L T$ is an isomorphism. Observe that these are torsors for the fpqc topology, but not for the fppf one: it may happen that, since $G$ is not of finite type, such a $T$ is not trivialized by any finite extension of $k$.
\end{itemize}
We have a natural map
\[\tau_L:\H^1_{\rm fpqc}(L,G)=\on{Tors}(L,G)\to\projlim_i\H^1(L,G_i)=\projlim_i\on{Tors}(L,G_i).\]
We are going to prove that, if $G$ is pro-finite, the map $\tau_L$ is an isomorphism. 

For completeness, we also remark that if $G$ is abelian these cohomologies coincide with fppf continuous cohomology in the sense of Jannsen. The proof follows directly from Jannsen's definition.

\begin{lemma}\label{schimg}
	Let $G=\projlim_iG_i$ a pro-finite group scheme over a field $k$, and let $G_i'$ be the scheme theoretic image of $G\to G_i$. Then $G_j\to G_i$ factorizes as $G_j\to G_i'\to G_i$ for every $j\gg i$ great enough.
	\begin{proof}
		If $G_i$ was finite étale, this would reduce to the analogous fact for finite groups, which in turn follows from the fact that a projective system of finite, non-empty sets is non-empty.
		
		For finite group schemes, write $G_i=\spec A_i$, then $G$ is the spectrum of
		\[A=\injlim_i A_i=\bigsqcup_iA_i/\sim\]
		where $\sim$ is the equivalence relation that identifies an element $a\in A_i$ with its image in $A_j$ for every $i\le j$. Call $K_{ij}\s A_i$ the kernel of $A_i\to A_j$ for every $i\le j$, $K_{ij}$ increases with $j$ and since $A_i$ is a finite dimensional vector space over $k$ we have that $K_{ij}$ is eventually stable, call $K_i$ the stable kernel. It is immediate to check that $K_i$ is the kernel of $A_i\to \bigsqcup_iA_i/\sim$, thus $G_i'=\spec A_i/K_i$ and the claim follows.
	\end{proof}
\end{lemma}

\begin{lemma}\label{prosurj}
	Let $G=\projlim_iG_i$ a pro-finite group scheme over a field $k$, and let $G_i'$ be the scheme theoretic image of $G\to G_i$, $(G_i')_i$ defines a second projective system of finite group schemes. It is possible to define a third projective system $(G_i)_i\sqcup(G_i')_i$ such that $(G_i)_i$ and $(G_i')_i$ are both cofinal sub-systems.
	\begin{proof}
		Let $I$ be the poset of indexes of $(G_i)_i$, and consider a copy of it $I'=I$. On the disjoint union $I\sqcup I'$, define an order in the following way.
		
		The restriction of the order to each component is just the order on $I$. If $i\in I$ and $j'\in I'$ (corresponding to $j\in I$), then $j'\ge i$ if $j\ge i$. If $i'\in I'$ and $j\in I$, then $j\ge i'$ if $j\ge i$ and the morphism $G_j\to G_i$ factorizes as $G_j\to G_i'\to G_i$. It is obvious that $I'$ is cofinal in $I\sqcup I'$. For every $i$ and for every $j>>i$ great enough, we have that $G_j\to G_i$ factorizes as $G_j\to G_i'\to G_i$ thanks to \autoref{schimg}: this tells us that $I$ is cofinal in $I\sqcup I'$, too.  
	\end{proof}
\end{lemma}

The proof of the following is due to A. Vistoli.

\begin{proposition}\label{proftor}
	Let $G=\projlim_i G_i$ a pro-finite group. Then 
	\[\tau:\on{Tors}(L,G)\to\projlim_i\H^1(L,G_i)=\projlim_i\on{Tors}(L,G_i)\]
	is an isomorphism.
	\begin{proof}
		Thanks to \autoref{prosurj}, we may suppose that $G\to G_i$ is surjective for every $i$, by which we mean that the associated morphism of Hopf algebras is injective.
	
		First, let us prove surjectivity of $\tau$. Consider an element $(T_i)_i\in\projlim_i\on{Tors}(L,G_i)$, for every $i\le j$ let $H_{ij}$ be the set of $G_j\to G_i$-equivariant morphisms $\sigma_{ij}:T_j\to T_i$. By hypothesis $H_{ij}$ is nonempty, we want to show that we can choose one $\sigma_{ij}\in H_{ij}$ for every $i\le j$ such that $\sigma_{ij}\circ\sigma_{jk}=\sigma_{ik}$ for every $i\le j\le k$. Consider $H=\prod_{i\ge j}H_{ij}$, we have that $H_{ij}$ is finite and thus if we consider $H$ with the product topology, it is compact. For $a\le b\le c$, let $C_{abc}\s H$ be subset of $(\sigma_{ij})_{ij}\in H$ such that $\sigma_{ab}\circ\sigma_{bc}=\sigma_{ac}$, it is a closed subset. The claim is equivalent to showing that $\bigcap_{a\ge b\ge c}C_{abc}$ is non-empty: since $H$ is compact, it is enough to check that the intersection of a finite number of them is non-empty. 
		
		Let $S$ be a finite set of triplets $a\ge b\ge c$, and choose an index $l$ such that $l\ge a$ for every triplet in $S$. For every $i\le l$, choose any $\sigma_i:T_l\to T_i$, and for every $i\le j\le l$ define $\sigma_{ij}:T_j\to T_i$ as the only equivariant morphism such that $\sigma_{ij}\circ\sigma_j=\sigma_i$: the unicity follows from the fact that $G_l\to G_j$ is surjective. Then we have that
		\[\sigma_{ij}\circ\sigma_{jk}\circ\sigma_{k}=\sigma_{ij}\circ\sigma_{j}=\sigma_i,\]
		and thus $\sigma_{ij}\circ\sigma_{jk}=\sigma_{ik}$ as desired.
		
		For injectivity, let $T,T'$ two $G$-torsors such that $T_i\simeq T_i'$ for every $i$. Let $H_i$ be the set of $G_i$-equivariant isomorphisms $T_i\simeq T_i'$, $H_i$ is finite for every $i$. This makes $(H_i)$ into a projective system of finite, non-empty sets, thus its limit is non-empty, and this allows us to define an isomorphism $T\simeq T'$. 
	\end{proof}
\end{proposition}

\subsection{Closed subgroups and essential dimension in the pro-finite case}

If $H\s G$ are finite group schemes, it is known that $\ed H\le\ed G$. The usual proof of this uses a generically free representation of $H$, see \cite[Theorem 6.19]{bf03}, but we do not have such representations for pro-finite group schemes.

There exists another, more direct proof, see \cite[Proposition 2.17]{brv07}: this one works in the pro-finite case. Some of their hypotheses are not satisfied by pro-finite group schemes, though, so let us write the details.

\begin{lemma}\label{subgrpinfed}
	Let $G$ be a pro-finite group scheme over a field $k$ with a closed sub-group $H\s G$. Let $T\to\spec L$ be an $H$-torsor. Then $\ed T\times^{H}G=\ed T$.
	\begin{proof}
		Clearly, $\ed T\times^{H}G\le\ed T$, we want to prove that $\ed T\times^{H}G\ge\ed T$.
		
		Let $k\s E\s L$ be a subextension with a $G$-torsor $Q\to\spec E$ such that $Q_{L}\simeq T\times^{G}H$, we want to show that $\trdeg (E/k)\ge\ed T$. The liftings of $Q$ to $H$ are represented by the pro-finite étale scheme $P=\on{Mor}_{G}(Q,G/H)$ over $E$.
		
		The fact that $T$ lifts $Q_{L}$ to $H$ gives us a point $\spec L\to P$ with residue field $L_{0}\s L$, this corresponds to an $H$-torsor $T'\to\spec L_{0}$ with $T'_{L}=T$. In particular, $\trdeg (L_{0}/k)\ge\ed T$. Since $P$ is pro-finite over $E$, it follows that $L_{0}$ is algebraic over $E$, and thus $\trdeg (E/k)=\trdeg (L_{0}/k)\ge\ed T$.
	\end{proof}
\end{lemma}

\begin{corollary}
	Let $G$ be a pro-finite group scheme over a field $k$ with a closed sub-group $H\s G$. Then $\ed H\le\ed G$.\qed
\end{corollary}

\section{Essential dimension of pro-finite group schemes}\label{infsect}

In this section, we prove that a pro-finite étale group scheme $G$ has finite essential dimension if and only if it is finite, in a strong sense.

\begin{theorem}\label{crit}
	Let $G$ be a pro-finite, infinite étale group scheme over a field $k$. Then there exists a $G$-torsor of infinite essential dimension.
	
	In particular, a pro-finite étale group scheme is finite if and only if it has finite essential dimension.
\end{theorem}

If $G$ is an algebraic group, then it is known that $\ed G<\infty$, see for instance \cite[Proposition 4.11, Remark 4.12]{bf03}. We only need to prove the part of \autoref{crit} regarding infinite groups. 

The proof is divided in two parts: first we prove the theorem for a few basic abelian cases, then we show that every infinite, pro-finite group contains a copy of one of these and hence conclude by \autoref{subgrpinfed}.  

\subsection{The basic abelian cases}

The statement of \autoref{crit} for the basic abelian cases is implied by the following \autoref{killer}. The hypotheses of \autoref{killer} are very restrictive and could be relaxed with some additional work, but this is not necessary: \autoref{killer} is just an intermediate step towards \autoref{crit}, and the few groups for which we need it already satisfy its hypotheses.

\begin{lemma}\label{killer}
	Let $G=\projlim_{n\in\N}G_{n}$ be a pro-finite abelian group scheme over an infinite field $k$, and let $H_{n}=\ker(G_{n}\to G_{n-1})$. Suppose that for every extension $L/k$ and every $n$
	\[\H^{1}(L,G_{n+1})\to\H^{1}(L,G_{n})\]
	is surjective and that for every $n$ there exists an $H_{n}$-torsor on $k(t)$ of essential dimension $1$. Then there exists a $G$-torsor $T\to\spec k(t_{1},t_{2},\dots)$ such that $\ed_{k} T=\infty$.
	\begin{proof}
		We break the proof in two steps. For every $n$, write $k(\tau_{n})=k(t_{1},\dots,t_{n})$ and $k(\tau_{\infty})=k(t_{1},t_{2},\dots)$. Recall that $\H^{1}$, for abelian groups, has an abelian group structure, i.e. we can add and subtract torsors for abelian group schemes.
			
		{\bf Step 1. Definition of $T_{n}$.} For every $n$, let $S'_{n}\to\spec k(t_{n})$ be some $H_{n}$-torsor with $\ed_{k} S'_{n}=1$, set $S_{n}=S'_{n}\times^{H_{n}}G_{n}$. We have that $\ed S_{n}=1$ thanks to \autoref{subgrpinfed}.
						
		Now define the $G_{n}$-torsor $T_{n}\to\spec k(\tau_{n})$ recursively in the following way. For $n=1$, set
		\[T_{1}=S_{1}.\]
		For $n\ge 2$, choose $T'_{n}\to\spec k(\tau_{n-1})$ as any lifting of $T_{n-1}$ (which exists by hypothesis) and define
		\[T_{n}=T'_{n,k(\tau_{n})}+S_{n,k(\tau_{n})}.\]
		Attention: \emph{first} we lift $T_{n-1}$ to $T'_n$ on the field $k(\tau_{n-1})$, \emph{then} we extend to $k(\tau_n)$ and finally we add $S_{n,k(\tau_{n})}$. If we had first extended to $k(\tau_n)$ and then lifted to $G_n$, then $S_n$ would not have played any role, $T_n$ and $T'_{n,k(\tau_{n})}$ would have been indistinguishable liftings of $T_{n-1}$. This is not the case: $T'_{n,k(\tau_{n})}$ descends to $k(\tau_{n-1})$ while $T_n$ does not.     
		
		{\bf Claim 1:} $T_n$ does not descend to $k(\tau_{n-1})$.
		
		Since $T'_{n,k(\tau_{n})}$ descends to $k(\tau_{n-1})$, this is equivalent to showing that $S_{n,k(\tau_{n})}$ does not descends to $k(\tau_{n-1})$.
		
		By contradiction, suppose we have a $G_n$-torsor $R_n$ on $k(\tau_{n-1})$ which extends to $S_{n,k(\tau_{n})}$. Recall that $S_n$ is a $G_n$-torsor on $k(t_n)$ which extends to $S_{n,k(\tau_{n})}$.
		
		Let $\tilde{R}_n$ and $\tilde{S}_n$ be $G_{n}$-torsors which are spreading outs of $R_{n},S_{n}$ over open subsets $U,V$ of $\A^{n-1}$, $\A^{1}$. Up to shrinking $U,V$, we may suppose that we have an isomorphism
		\[U\times\tilde{S}_{n}\simeq\tilde{R}_{n}\times V\to U\times V\s\A^{n}\]
		since this is true generically. Since $k$ is infinite, we can choose a rational point $u\in U(k)$. If we restrict the equality above to $u\times\spec k(t_{n})\in U\times V$, we get
		\[S_{n}\simeq\tilde{R}_{n,u}\times_k\spec k(t_{n}).\]
		But $\ed_{k} S_{n}=1$ by construction, hence we have a contradiction. This concludes the first step.
		
		{\bf Step 2. Definition of $T$.} It is clear that the pushforward of $T_{n}$ to $G_{n-1}$ is $T_{n-1,k(\tau_{n})}$, since $S_{n,k(\tau_{n})}$ maps to $0$. Define $T$ as the projective limit 
		\[T=\projlim_n T_{n,k(\tau_{\infty})}\in\H^{1}(k(\tau_{\infty}),G)=\projlim_{n}\H^{1}(k(\tau_{\infty}),G_{n}).\] 
		
		{\bf Claim 2:} $\ed_k T=\infty$.
		
		Suppose by contradiction that $T$ descends to some sub-extension of finite transcendence degree of $k(\tau_{\infty})$. In particular, $T$ descends to $k(\tau_m)$ for some $m$ large enough, thus $T_{m+1,k(\tau_{\infty})}=T\times^{G}G_{m+1}$ descends to $k(\tau_m)$ too.
		
		Let $Q$ be a $G_{m+1}$-torsor over $k(\tau_m)$ which extends to $T_{m+1,k(\tau_{\infty})}$. Consider the difference $G_{m+1}$-torsor
		\[D=Q_{k(\tau_{m+1})}-T_{m+1}.\]
		By construction, $D_{k(\tau_{\infty})}$ is trivial, but since $k(\tau_{m+1})$ is algebraically closed in $k(\tau_{\infty})$ we have that $D$ itself is trivial, i.e.
		\[Q_{k(\tau_{m+1})}\simeq T_{m+1}.\]
		But this gives a contradiction, since $T_{m+1}$ does not descend to $k(\tau_m)$ by claim 1.
	\end{proof}
\end{lemma}

\begin{lemma}\label{mut}
	Over any field $k$ and for every prime $l$ (possibly equal to $\on{char}k$), there exists a $\mu_{l}$-torsor on $k(t)$ of essential dimension $1$.
	\begin{proof}
		Observe that if the characteristic of $k$ is $l$, then the étale/Galois cohomology of $\mu_{l}$ is trivial. However, the set of $\mu_{l}$-torsors corresponds to fppf cohomology, and fppf cohomology can be different from étale cohomology for non-smooth groups.
		
		For fppf cohomology, regardless of the characteristic of $k$, elementary Kummer theory provides us with a natural isomorphism
		\[\H^{1}(k(t),\mu_{l})=k(t)^{*}/k(t)^{*l},\]
		and $t\in k(t)^{*}/k(t)^{*l}$ has essential dimension $1$ since it does not descend to $k$, which is the only subextension of transcendence degree $0$. This can be checked using the valuation $v:k(t)^{*}\to\Z$ that sends $t$ to $1$: in fact, $v$ induces an homomorphism $k(t)^{*}/k(t)^{*l}\to\Z/l$ trivial on $k$ but not trivial on $t$.
	\end{proof}
\end{lemma}

\begin{corollary}\label{edzp1}
	Over any field $k$ and for any prime $l$, there exists a $\Z_{l}(1)$-torsor $T\to\spec k(t_{1},t_{2},\dots)$ such that $\ed_{k} T=\infty$.
	\begin{proof}
		Write $\Z_{l}(1)=\projlim_{n\in\N}\mu_{l^{n}}$. The map $\H^1(L,\mu_{l^{n+1}})\to\H^1(L,\mu_{l^{n}})$ is surjective by Kummer theory for every $L/k$ and every $n$: along with \autoref{mut}, this ensures that the hypotheses of \autoref{killer} are satisfied.
	\end{proof}
\end{corollary}

\begin{lemma}\label{pnt}
	If $\on{char}k=p>0$, there exists a $\Z/p$-torsor on $k(t)$ of essential dimension $1$.
	\begin{proof}
		Let $\Phi:k(t)\to k(t)$ the homomorphism  $x\mapsto x^{p}-x$, using the Artin-Schreier exact sequence we have
		\[\H^{1}(k(t),\Z/p)\simeq k(t)/\Phi(k(t)).\]
		Let us prove that $t\in k(t)/\Phi(k(t))$ has essential dimension $1$. In fact, if $t$ is defined on $k$ (the only algebraic sub-extension of $k(t)/k$) we have
		\[t=\lambda+q^{l}-q\]
		for some $\lambda\in k$ and $q\in k(t)$. But then $k(t)\s k(q)\s k(t)$ and thus $q(t)=(at+b)/(ct+d)$ is a non-constant linear transformation, we get the equation
		\[-c^pt^{p+1}+(ac^{p-1}-\lambda c^p-a^p)t^{p} + r(t)=0\]
		for some $r(t)\in k[t]$ of degree at most $p-1$. This implies that $c=a=0$ and hence $q=b/d$ is constant, which is absurd.
	\end{proof}
\end{lemma}

\begin{corollary}\label{edzp}
	If $k$ is a field of characteristic $p>0$, there exists a $\Z_{p}$-torsor $T\to\spec k(t_{1},t_{2},\dots)$ such that $\ed_{k} T=\infty$.
	\begin{proof}
		Let us check the hypotheses of \autoref{killer}. The surjectivity of $\H^{1}(L,\Z/p^{n+1})\to\H^{1}(L,\Z/p^{n})$ follows from the fact that $\H^{2}(L,\Z/p)$ is trivial, see \cite[Ch. 2, Proposition 3]{ser94}. The kernel of $\Z/p^{n+1}\to\Z/p^{n}$ is just $\Z/p$, thanks to \autoref{pnt} there exists a $\Z/p$-torsor on $k(t)$ of essential dimension $1$.
	\end{proof}
\end{corollary}

\begin{corollary}\label{infp}
	Let $k$ be a field, $l_1,l_2,\dots$ a sequence of not necessarily distinct prime numbers. Let
	\[G=\prod_{i=1}^{\infty}\Z/l_{i}.\]
	Then there exists a $G$-torsor of infinite essential dimension.
	\begin{proof}
		We may extend $k$ and suppose $k=\bar{k}$, and hence $\mu_{l}=\Z/l$ if $\on{char}k\neq l$. Write
		\[\prod_{i=1}^{\infty}\Z/l_{i}=\projlim_{n\in\N}\prod_{i=1}^{n}\Z/l_{i},\]
		and let us check the hypotheses of \autoref{killer}.
		
		The surjectivity of
		\[\H^{1}\left(L,\prod_{i=1}^{n}\Z/l_{i}\right)\to\H^{1}\left(L,\prod_{i=1}^{n-1}\Z/l_{i}\right)\]
		is obvious, since cohomology commutes with direct product. The fact that there exists a $\Z/l_{i}$-torsor of essential dimension $1$ on $k(t)$ comes either from \autoref{mut} if $\on{char} k\neq l_{i}$ or from \autoref{pnt} if $\on{char} k=l_{i}$.
	\end{proof}
\end{corollary}

\subsection{The general case}

\begin{lemma}\label{infabsub}
	Let $G$ be an infinite pro-finite group. At least one of the following is true.
	\begin{itemize}
		\item There exists a closed sub-group isomorphic to $\Z_{l}$ for some prime $l$.
		\item There exists a closed sub-group isomorphic to $\prod_{i\in\N}\Z/l_{i}$ for some primes $l_{i}$.
	\end{itemize}
	\begin{proof}
		By a theorem of Zel'manov \cite[Theorem 2]{ze92}, $G$ contains an infinite abelian subgroup. Hence, we may suppose that $G$ is abelian.
		
		If $G$ is abelian, then $G=\prod_{l}G(l)$ is the product of its $l$-Sylow subgroups. If all the $G(l)$ are finite, then it is easy to construct a subgroup of the form $\prod_{i\in\N}\Z/l_{i}$ by choosing a cyclic subgroup of prime order of $G(l)$ for every $l$. Otherwise, there exists an $l$ with $G(l)$ infinite, we may replace $G$ with $G(l)$ and assume that $G$ is pro-$l$.
		
		If $G$ is pro-$l$ and contains an element which is not torsion, then we have a closed subgroup $\Z_{l}\s G$.

		Suppose now that $G$ is pro-$l$ and torsion. First, let us show that $G[l]$ is infinite. If $G[l]$ is finite, then by induction on the multiplication by $l$ map $G[l^{n+1}]\to G[l^{n}]$ we get that $G[l^{n}]$ is finite for every $n$. Since a pro-$l$ group has no non-trivial $l$-divisible element, there exists an $n_{0}$ such that the multiplication by $l^{n}$ map $G[l^{n+1}]\to G[l]$ is trivial for every $n\ge n_{0}$. It follows that $G[l^{n}]=G[l^{n_{0}}]$ for every $n\ge n_{0}$, hence $G=G[l^{n_{0}}]$, which gives a contradiction since $G$ is infinite.
		
		Thus, $G$ is a torsion, abelian pro-$l$ group with $G[l]$ infinite. Let $H$ be the Pontryagin dual of $G[l]$, it is an $l$-torsion abelian group. By \cite[Theorem 6]{kap54}, we may write $H\simeq\bigoplus_{I}\Z/l$ for some infinite set $I$. Let $S\s I$ be an infinite, countable subset, we have a quotient map $H\to\bigoplus_{S}\Z/l$. The Pontryagin dual $\prod_{S}\Z/l$ of $\bigoplus_{S}\Z/l$ is a closed subgroup of $G[l]$.
	\end{proof}
\end{lemma}

Let us now prove \autoref{crit}.

\begin{proof}[Proof of {\autoref{crit}}]
	We may extend $k$ and suppose that $k=\bar{k}$, hence we may regard $G$ simply as a pro-finite group. By \autoref{infabsub}, there exists an infinite closed abelian subgroup $H\s G$ which is isomorphic either to $\Z_{l}$ for some $l$ or to an infinite product of the form $\prod_{i\in\N}\Z/l_{i}$. Either by \autoref{edzp1}, \autoref{edzp} or \autoref{infp}, there exists an $H$-torsor of infinite essential dimension. The statement then follows by \autoref{subgrpinfed}.
\end{proof}

In view of \autoref{crit}, essential dimension is not a very interesting invariant for infinite, pro-finite étale group schemes. One may be content with this, and be done with it. However, we think that the ideas and formalism of essential dimension may still give non-trivial information, at the cost of modifying the basic definition of essential dimension.

\section{Finite type essential dimension}

The results of \autoref{infsect} are based on \autoref{killer}, where we construct a single torsor of infinite essential dimension. Observe that the proof of \autoref{killer} does \emph{not} adapt to the construction of torsors with finite and arbitrarily large essential dimension: we really use the "gap" between finite and infinite. Moreover, in order to apply \autoref{killer} in the proof of \autoref{crit}, we need to pass to the algebraic closure of the base field, thus "killing" the Galois structure of $G(\bar{k})$.

We define finite type essential dimension by focusing only on finitely generated extension. This has the advantages of avoiding the pathology of \autoref{killer} and preserving the arithmetic structure of the group.

\begin{definition}
	Let ${\rm F}:\fld_{k}\to {\rm Set}$ be a functor from the category of extensions of $k$ to ${\rm Set}$. The \emph{finite type essential dimension} $\fed_{k}\rm{F}$ is the supremum of the essential dimensions $\ed_{k}(\alpha)$ where $\alpha$ varies among objects $\alpha\in {\rm F}(L)$ with $L$ a finitely generated extension of $k$.
\end{definition}

Note that in the original definition of $\ed G$ for groups $G$ of finite type in \cite{rei00} only $G$-torsors over finitely generated fields $K/k$ were considered, since the point view was the one of group actions on varieties. However, the restriction on $K$ was later dropped, since for groups of finite type $\ed G=\fed G$.

\begin{remark}
	To remain in Merkurjev's general framework of essential dimension for functors $\fld_{k}\to {\rm Set}$ (see \cite{bf03}) one can give the following alternative definition, for which I thank Z. Reichstein.
	
	Given a functor ${\rm F}:\fld_{k}\to {\rm Set}$, we may define the functor $\rm{F^{fin}}:\fld_{k}\to {\rm Set}$ as
	\[{\rm F^{fin}}(L)=	\begin{cases}
							{\rm F}(L)	&	L/k\text{ is finitely generated}	\\
							\{\ast\}		&	\text{otherwise}
						\end{cases}\]
	Then we have $\fed_k{\rm F}=\ed_k{\rm F^{fin}}$.
\end{remark}

\begin{lemma}
	Let $G$ be a pro-finite group scheme over a field $k$ with a closed sub-group $H\s G$. Then $\fed H\le\fed G$.
	\begin{proof}
		This follows from \autoref{subgrpinfed}.
	\end{proof}
\end{lemma}

\begin{lemma}
	Let $G$ be a group scheme of finite type over $k$, then $\fed_k G=\ed_k G$.
	\begin{proof}
		This follows from the fact that every $G$-torsor descends to a field finitely generated over $k$.
	\end{proof}
\end{lemma}

\begin{example}
	Let $X$ be a smooth, projective curve over a field $k$ finitely generated over $\Q$. Suppose it has a rational point $x\in X(k)$ and let $\upi_{1}(X,x)$ be its étale fundamental group scheme, it is infinite and thus $\ed_{k}\upi_{1}(X,x)=\infty$.
	
	On the other hand, if Grothendieck's section conjecture is true then 
	\[X(k')=\H^{1}(k',\upi_{1}(X,x))\]
	for every finitely generated extension $k'/k$, and thus
	\[\fed_{k}\upi_{1}(X,x)=\fed_{k} X=1.\]
\end{example}

\begin{remark}
	Excluding the category of fields finitely generated over the base field, there is another category of extensions of the base field that might have been a good candidate, i.e. the category of fields of finite transcendence degree over the base field. However, using this category would not solve our problem: the modified essential dimension relative to this category would still be infinite too often, see \autoref{faked}.	
\end{remark}

Recall that, if $l$ is a prime and $G$ is a pro-$l$ group scheme over a field $k$ with $\on{char}k\neq l$, we say that $G$ has a geometric subgroup if there exists a finite extension $k'/k$ and a semi-abelian variety $A$ over $k'$ with a non-trivial homomorphism $T_{l}A\to G_{k'}$.

We are going to prove the following.

\begin{theorem}\label{fed}
	Let $l$ be a prime and $G$ an abelian, torsion free, pro-$l$ group scheme over a field $k$ with $\on{char}k\neq l$. If $G$ has a geometric subgroup, then $\fed G=\infty$, otherwise $\fed G=0$.
\end{theorem}

Both cases of \autoref{fed} are non-trivial. We thus split \autoref{fed} in the two following theorems \ref{fedsab} and \ref{fed0}.

\subsection{Groups with infinite $\fed$}

\begin{theorem}\label{fedsab}
	Let $G$ be an abelian, torsion free pro-$l$ group scheme over a field $k$ with $\on{char}k\neq l$. If $G$ has a geometric subgroup, then $\fed_{k}G=\infty$.
	\begin{proof}
		Since $G$ has a geometric subgroup, there exists a finite extension $k'$ and a semi-abelian variety $A$ over $k'$ with a non-trivial homomorphism $\omega:T_{l}A\to G_{k'}$. Up to passing to a finite extension of $k$ and to subgroups of $A$, we may suppose that $k=k'$ and that $A$ is either $\mathbb{G}_{m}$ or a geometrically simple abelian variety. Let $g$ be the dimension of $A$, and fix any positive integer $d$, we are going to construct a $G$-torsor of essential dimension $dg$. In order to deal better with some base-point issues, we will use the language of gerbes.
		
		Let $D=\hom(A,A)\otimes_{\Z}\Q$, it is a division algebra finite over $\Q$. Composition with $\omega$ gives a $\Q_{l}$-linear map
		\[D_{\Q_{l}}\to\hom(T_{l}A,T_{l}A)_{\Q_{l}}\to\hom(T_{l}A,G)_{\Q_{l}}\to\H^{1}(A,G)_{\Q_{l}}.\]
		Since $D$ has finite dimension over $\Q$, there exists a finitely generated subextension $\Q\s L\s \Q_{l}$ and an $L$-vector space $W\s D_{L}$ whose extension $W_{\Q_{l}}$ is the kernel of $D_{\Q_{l}}\to\H^{1}(A,G)_{\Q_{l}}$.
		
		Since $\Q_{l}$ is not finitely generated over $\Q$, then $[\Q_{l}:L]=\infty$ and we may choose a $d$-uple of $l$-adic integers $\underline{n}=(n_1,\dots,n_d)\in(\Z_{l}^*)^{\oplus d}$ which are linearly independent over $L$. Define an homomorphism
		\[\underline{n}\omega:T_{l}A^{d}=\oplus_{d}T_{l}A\to G\]
		\[(a_{1},\dots,a_{d})\mapsto n_{1}\omega(a_{1})+\dots+n_{d}\omega(a_{d}).\]	
		The homomorphism $\underline{n}\omega$ defines a $G$-torsor on $A^{d}$ which we still call $\underline{n}\omega$ with an abuse of notation. We want to show that $\ed_{k}\underline{n}\omega\ge dg$, since $\ed_{k}\underline{n}\omega\ge \ed_{\bar{k}}\underline{n}\omega_{\bar{k}}$ we may assume $k=\bar{k}$.
		
		Suppose by contradiction that $\underline{n}\omega$ descends to a subextension $k(A^d)/k'/k$ of strictly smaller transcendence degree. Since $k=\bar{k}$ is perfect, by generic smoothness we may find a smooth variety $V$ with $k(V)=k'$ and an open subset $U\s A^{d}$ with a surjective morphism $U\to V$ extending $\spec k(A^{d})\to\spec k'$. By \autoref{ramext}, the morphism $\spec k'\to BG$ extends to $V\to BG$.
		
		Let $C\s U$ be any irreducible curve that gets contracted by $U\to V$, fix $\bar{C}$ a normalization, we have a natural morphism $f:\bar{C}\to A^{d}$ and the composition $\bar{C}\to C\to A^{d}\xar{\underline{n}\omega}BG$ factorizes through a rational section. Let $c\in \bar{C}(k)$ be any rational point, up to translating the morphism $\bar{C}\to A^{d}$ we may assume that $c$ maps to the origin $a\in A^{d}$. The translated morphism still gets contracted by $A^{d}\to BG$, since $A^{d}(\bar{C})\to \H^{1}(\bar{C},G)$ is a group homomorphism and translation by a point of $\H^{1}(k,G)$ sends the subgroup $\H^{1}(k,G)\s \H^{1}(\bar{C},G)$ to itself.
		
		Let $\bar{C}\to J$ be the semi-abelian Jacobian of $\bar{C}$ and $T\s J$ the toric part, by construction $T_{l}J\to G$ is $0$. At least one coordinate of $J\to A^{d}$ is non-trivial, say the first one, let $e_{1}:A^{d}\to A$ be the first projection. If $A=\mathbb{G}_{m}$ then $T\to A^d\xar{e_{1}} A$ is non-trivial, too, because there are no non-trivial homomorphisms $J/T\to A$. If $A$ is an abelian variety, then we may replace $\bar{C}$ with its completion (since a generically constant étale torsor is constant) and assume that $T=\{0\}$. In both cases, there exists an homomorphism $A\to J$ such that the composition $A\to J\to A^{d}\xar{e_{1}} A$ is an isogeny.
				
		Call $\underline{\phi}=(\phi_{1},\dots,\phi_{d}):A\to A^{d}$ the composition, $\phi_{1}$ is an isogeny and the composition
		\[\underline{n}\omega\circ T_{l}\underline{\phi}:T_{l}A\to T_{l}A^{d}\to G\]
		is $0$.
		
		We have that the image $[\phi_{1}]$ of $\phi_{1}\neq 0\in D$ in $D_{\Q_{l}}/W_{\Q_{l}}\s\H^{1}(A,G)_{\Q_{l}}$ is non-zero: it is associated to $\omega\circ T_{l}\phi_{1}:T_{l}A\to T_{l}A\to G$, which is non-trivial since $\omega:T_{l}A\to G$ is non-trivial, $T_{l}\phi_{1}$ has maximal rank ($A$ is simple) and $G$ is torsion free.
		
		Recall that we have a finitely generated subextension $\Q\s L\s\Q_{l}$, an $L$-linear subspace $W\s D_{L}$ and an injective homomorphism
		\[D_{\Q_{l}}/W_{\Q_{l}}=(D_{L}/W)\otimes_{L}\Q_{l}\to\H^{1}(A,G)_{\Q_{l}}.\]
		Moreover, $n_{1},\dots,n_{d}\in\Q_{l}$ are linearly independent over $L$, $[\phi_{i}]\in D_{L}/W$ for $i=1,\dots,n$ and $[\phi_{1}]\neq 0$. It follows that the linear combination
		\[[\phi_{1}]\otimes n_{1}+\dots+[\phi_{d}]\otimes n_{d}\in (D_{L}/W)\otimes_{L}\Q_{l}\s\H^{1}(A,G)_{\Q_{l}}\]
		is non-zero. But this is absurd, since this linear combination is precisely $\underline{n}\omega\circ T_{l}\underline{\phi}$, which is $0$ by construction.
		
	\end{proof}
\end{theorem}

Using \autoref{fedsab}, we can show that using the category of field extensions of finite transcendence degree still leads to a trivial notion of essential dimension.

\begin{corollary}\label{faked}
	Let $G$ be a pro-finite étale group scheme. Suppose that there exists a prime $l\neq\on{char}k$ such that an $l$-Sylow subgroup of $G(\bar{k})$ is infinite. Then there exist $G$-torsors of finite but arbitrarily large essential dimension.
	\begin{proof}
		By \autoref{infabsub}, there exists a closed subgroup $G'\s G(\bar{k})$ which isomorphic to either $\Z_{l}$ or an infinite product $\prod \Z/l$. Passing to $\bar{k}$ and using \autoref{subgrpinfed} to replace $G$ with $G'$, we may suppose $k=\bar{k}$ and that $G$ is either $\Z_{l}$ or an infinite product $\prod_{\N}\Z/l$.
		
		If $G=\Z_{l}$, since $\on{char}k\neq l$ and $k=\bar{k}$, $\Z_{l}\simeq\Z_{l}(1)$ and thus we obtain the statement by \autoref{fedsab}.
		
		If $G$ is an infinite product $\prod_{\N}\Z/l$, then for every $n$ we have a closed subgroup $\prod_{n}\Z/l\s G$. Since $\on{char} k\neq l$ and $k=\bar{k}$, we have $\ed\prod_{n}\Z/l=n$ and thus we conclude by \autoref{subgrpinfed}.
	\end{proof}
\end{corollary}

\subsection{Groups with $\fed=0$}

Let us now address the case in which $G$ has no geometric subgroups. First, we need to prove that if $G$ has no geometric subgroups, this remains true after a finitely generated extension of the base field.

Recall that a field extension $k'/k$ is regular if $k$ is algebraically closed in $k'$ and there exists an intermediate extension $k'/h/k$ such that $h/k$ is purely transcendental and $k'/h$ is separable algebraic. If $k'/k$ is a regular extension and $A$ is an abelian variety over $k'$, the $k'/k$-\emph{trace} of $A$ is a final object in the category of abelian varieties $B$ over $k$ with an homomorphism $B_{k'}\to A$: the $k'/k$-trace always exists, see \cite[Theorem 6.2]{co06}. If furthermore $k'/k$ is finitely generated and $B$ is the $k'/k$-trace, the Lang-Néron theorem says that $A(k')/B(k)$ is finitely generated, see \cite[Theorem 7.1]{co06}. 

\begin{lemma}\label{tatequot}
	Let $k'/k$ be a finitely generated regular extension, $A$ a geometrically simple abelian variety over $k'$, $l$ a prime different from $\on{char}k$. The following are equivalent.
	\begin{enumerate}[(i)]
		\item $A$ descends to $k$.
		\item There exists a torsion free pro-$l$ group scheme $G$ over $k$ with a non-trivial homomorphism $T_{l}A\to G_{k'}$.
		\item There exists a non-trivial closed subgroup $H\s T_{l}A$ which descends to $k$.
	\end{enumerate}
	\begin{proof}
		$(i)\Rightarrow (ii).$ If $A$ descends to an abelian variety $B$ over $k$, then we may take $G=T_{l}B$.
		
		$(ii)\Rightarrow (iii).$ Denote by $\gal(k),\gal(k')$ the absolute Galois groups, since $k$ is algebraically closed in $k'$ we have a surjective homomorphism $\gal(k')\to\gal(k)$. Let $G'\s G_{k'}$ be the image of $T_{l}A\to G_{k'}$, since the action of $\gal(k')$ restricts to the geometric points of $G'$ and $\gal(k')\to\gal(k)$ is surjective it follows that $G'$ descends to a closed subgroup of $G$, i.e. we may suppose that $T_{l}A\to G_{k'}$ is surjective. In particular, $G(\bar{k})$ is topologically finitely generated. 
		
		Denote by $\hat{G}=\hom(G,\Z_{l}(1))$, then we have an injective homomorphism $\hat{G}_{k'}\to \hat{T_{l}A}=T_{l}\hat{A}$. Any isogeny $\hat{A}\to A$ gives us an embedding $\hat{G}_{k'}\to T_{l}\hat{A}\to T_{l}A$.
		
		$(iii)\Rightarrow (i).$ Suppose by contradiction that $A$ is not defined over $k$, since it is simple in particular its $k'/k$-trace is trivial. If $h/k$ is any extension and $hk'$ is the fraction field of $h\otimes_{k} k'$ (it is a domain since $k'/k$ is regular), then the $h'/h$-trace of $A_{h'}$ is trivial too, see \cite[Theorem 6.6]{co06}. By the Lang-Néron theorem, $A(\bar{k}k')$ is finitely generated, and thus $T_{l}A(\bar{k}k')$ is trivial. On the other hand, $H(\bar{k}k')$ is non-trivial, since $H$ is non-trivial and descends to $k$. This gives a contradiction.
	\end{proof}
\end{lemma}

\begin{lemma}\label{regext}
	Let $k'/k$ be a finitely generated extension of fields, and assume that $k$ is separably closed in $k'$. There exist finite, purely inseparable extensions $h'/k'$, $h/k$ such that $h\s h'$ and $h'$ is regular over $h$.
	\begin{proof}
		Let $K$ be a perfect closure of $k$, then $K\otimes_{k}k'$ has one minimal prime and thus it makes sense to consider its fraction field $K'$, which is finitely generated over $K$ and purely inseparable over $k'$. Since $K$ is perfect, there exists a separating basis $x_{1},\dots,x_{n}\in K'$. Since $K'$ is algebraic over $k'$, there exists a finite subextension $k'\s K'_{0}\s K'$ such that $x_{1},\dots,x_{n}\in K'_{0}$. Let $K_{0}\s K'_{0}$ be the algebraic closure of $k$ in $K'_{0}$.
		
		Choose $S\s K_{0}'$ a finite set of generators for the extension $K_{0}'/K_{0}(x_{1},\dots,x_{n})$. Since $K'$ is separable over $K(x_{1},\dots,x_{n})$, the elements of $S$ are separable over $K(x_{1},\dots,x_{n})$ and hence there exists a finite, purely inseparable subextension $k\s h\s K$ such that the elements of $S$ are separable over $h(x_{1},\dots,x_{n})$. Define $h'=hK_{0}'=h(x_{1},\dots,x_{n},S)$: since $h'$ is generated over $h(x_{1},\dots,x_{n})$ by a set of separable elements, we have that $x_{1},\dots,x_{n}$ is a separable basis for $h'/h$. Moreover, $h'$ is generated over $k'$ by $hk'$ and $K_{0}'$, and these are both finite, purely inseparable extensions of $k'$, thus $h'/k'$ is finite and purely inseparable too.
	\end{proof}
\end{lemma}

\begin{corollary}\label{ffgquot}
	Let $G$ be a pro-$l$ group scheme over a field $k$ of characteristic different from $l$. The following are equivalent.
	\begin{enumerate}[(i)]
		\item There exists a finite extension $k'/k$ and a semi-abelian variety $A$ over $k'$ with a non-trivial homomorphism $T_{l}A\to G_{k'}$.
		\item There exists a finitely generated extension $k'/k$ and a semi-abelian variety $A$ over $k'$ with a non-trivial homomorphism $T_{l}A\to G_{k'}$.
	\end{enumerate}
	In particular, if $k'/k$ is finitely generated and $G$ has no geometric subgroups, then $G_{k'}$ has no geometric subgroups.
	\begin{proof}
		The implication $(i)\Rightarrow(ii)$ is obvious.
		
		Let $k'/k$ be finitely generated and $A$ a semi-abelian variety with a non-trivial homomorphism $T_{l}A\to G_{k'}$. By replacing $k$ with $\bar{k}^{k'}$, we may suppose that $k$ is algebraically closed in $k'$.
		
		Let $D\s A$ be the toric part. If $T_{l}D\to G_{k'}$ is non-trivial, up to a finite extension we may find a subgroup $\G_{m}\s D$ such that $\Z_{l}(1)=T_{l}\G_{m}\to G_{k'}$ is non-trivial. Since both $\Z_{l}(1)$ and $G$ are pro-finite étale and defined over $k$,  and $k$ is algebraically closed in $k'$, the homomorphism descends to $k$.
		
		If $T_{l}D\to G_{k'}$ is trivial, we have an induced homomorphism $T_{l}(A/D)\to G_{k'}$, i.e. we may suppose that $A$ is abelian. Up to a finite extension of $k'$, we may pass to an isogeny factor and suppose that $A$ is geometrically simple. Moreover, since $k$ is algebraically closed in $k'$ the image of $T_{l}A\to G_{k'}$ descends to a subgroup $G'$ of $G$, hence we may replace $G$ with $G'$ and suppose that $T_{l}A\to G_{k'}$ is surjective. In particular, $G(\bar{k})$ it topologically finitely generated, and its torsion subgroup is finite. If the torsion subgroup is not trivial, since it is finite it is easy to find a finite extension $k''/k$ and a non-trivial morphism $T_{l}\G_{m}=\Z_{l}(1)\to G_{k''}$.
		
		Suppose now that $G$ is torsion free. Thanks to \autoref{regext}, we may replace both $k$ and $k'$ with finite extensions and assume that $k'/k$ is regular. By \autoref{tatequot}, $A$ descends to an abelian variety $B$ over $k$. Since both $T_{l}B$ and $G$ are pro-finite étale and $k$ is algebraically closed in $k'$, the homomorphism $T_{l}A=T_{l}B_{k'}\to G_{k'}$ descends to an homomorphism $T_{l}B\to G$.
	\end{proof}
\end{corollary}

\begin{lemma}\label{insdes}
	Let $k'/k$ be a purely inseparable extension and $G$ a pro-finite étale group scheme over $k$. Then $\H^{1}(k,G)\to\H^{1}(k',G')$ is an isomorphism.
	\begin{proof}
		Since $G$ is étale, we may use Galois cohomology. The statement then follows from the fact that the homomorphism $\gal(k')\to\gal(k)$ of absolute Galois groups is an isomorphism.
	\end{proof}
\end{lemma}

\begin{theorem}\label{fed0}
	Let $G$ be an abelian pro-$l$ group scheme over a field $k$ with $\on{char}k\neq l$. If $G$ does not have a geometric subgroup, then $\fed_{k} G=0$.
	\begin{proof}
		We are going to use the language of gerbes. Let $K/k$ be a finitely generated extension and $\spec K\to BG$ a section, we want to show that it descends to $\bar{k}^{K}$, we may replace $k$ with $\bar{k}^{K}$ and assume that $k$ is algebraically closed in $K$. By induction, it is enough to prove the case $\trdeg(K/k)=1$. Thanks to \autoref{regext} and \autoref{insdes}, we may moreover suppose that $K/k$ is regular.
		
		Since $K/k$ is regular of transcendence degree $1$, there exists a smooth projective curve $X/k$ with function field $k(X)=k$. Thanks to \autoref{ffgquot} and \autoref{r-proper}, the generic section $\spec k(X)\to BG$ extends to a global section $X\to BG$. Let $X\to J$ be the Albanese torsor of $X$, since $G$ is abelian and $\on{char}k\neq l$, we have that $X\to BG$ extends to $J\to BG$. Let $k'/k$ be a finite separable extension such that $J(k')\neq\emptyset$, we may regard $J_{k'}$ as an abelian variety. We have that the induced homomorphism $T_{l}J_{k'}\to G_{k'}$ is trivial since $G$ has no geometric subgroups and thus $J_{k'}\to BG_{k'}$ is constant. But then $J\to BG$ is constant, too, thanks to \autoref{constbc}. The statement follows.
	\end{proof}
\end{theorem}

If $k$ is finitely generated over $\Q$, Deligne's theory of weights gives us a way of checking the hypothesis of \autoref{fed0}.

\begin{corollary}\label{fedwp}
	Let $k$ be a field finitely generated over $\Q$, and let $l$ be a prime number. Let $G$ be an abelian, torsion free pro-finite $l$-group scheme over $k$. Assume that the Galois representation $G\otimes_{\Z_{l}}\Q_{l}$ is pure of weight $w\neq-2,-1$. Then $\fed_{k} G=0$.\qed
\end{corollary}

If $k$ has all $l$-adic roots of the unity, then $\Z_{l}(1)\simeq\Z_{l}$ and thus $\fed_{k}\Z_{l}=\fed_{k}\Z_{l}(1)=\infty$. If $k$ is finitely generated over $\Q$, then it is an easy consequence of a theorem of Florence \cite[Theorem 4.1]{flo08} that
\[\lim_n\ed_{k}\Z/l^{n}\Z=\infty.\]

In view of this, the following is rather surprising.

\begin{corollary}\label{fedzp}
	Let $k$ be a field finitely generated over $\Q$, and let $l$ be a prime number. Then
	\[\fed_{k}\Z_{l}(n)=0\]
	for every $n\neq 1$.\qed
\end{corollary}

From \autoref{fedzp} for $n=0$, we obtain the fact that every $\Z_{l}$-extension of a field finitely generated over $\Q$ descends to a number field.

\begin{corollary}\label{galzp}
	Let $K$ be finitely generated over $\Q$, and let $k=\overline{\Q}^{K}$ be the algebraic closure of $\Q$ in $K$. If $H/K$ is a $\Z_{l}$-extension, there exists a $\Z_{l}$-extension $h/k$ such that $H=hK$. In particular, all $\Z_{l}$-extension of fields finitely generated over $\Q$ descend to number fields. 
	\begin{proof}
		We have that $\spec H\to\spec K$ is a $\Z_{l}$-torsor, by \autoref{fedzp} we have a $\Z_{l}$-torsor $\spec h\to\spec k$ such that $\spec H=\spec h\times_{k}\spec K$. In particular, $\spec h\to\spec k$ is connected and pro-finite étale, thus $h/k$ is $\Z_{l}$-Galois extension. The isomorphism $h\otimes_{k} K\simeq H$ allows us to fix an embedding $h\s H$ such that $hK=H$.
	\end{proof}
\end{corollary}

\subsection{Characteristic $p$}

Since most of our work about finite type essential dimension focuses on pro-$l$ groups in characteristic different from $l$, let us make here a brief remark about characteristic $p$. Let $\on{char}k=p>0$.

It can be easily shown that 
\[\ed_{k}\Z/p^{n}\le\ed_{k}\Z/p^{n+1}\le\ed_{k}\Z/p^{n}+1\]
using the obvious short exact sequence. A. Ledet, who obtained in \cite{le04} and \cite{le04b} the main known results about essential dimension of cyclic $p$-groups in characteristic $p$, made the following conjecture.

\begin{conjecture*}[Ledet]
	If $\on{char} k=p$, then $\ed_{k}\Z/p^{n}=n$, or equivalently $\ed_{k}\Z/p^{n}=\ed_{k}\Z/p^{n-1}+1$.
\end{conjecture*}

\begin{proposition}
	If $\on{char} k=p$ and $\lim_n\ed_{k}\Z/p^{n}=\infty$, then $\fed_{k}\Z_{p}=\infty$.
	\begin{proof}
		Since $\H^{2}(L,\Z/p)=0$ for every extension $L/k$, every $\Z/p^{n}$ torsor lifts to a $\Z/p^{n+1}$ torsor, and by recursion to a $\Z_{p}$-torsor. If $\lim_{n\to\infty}\ed_{k}\Z/p^{n}=\infty$, we have thus $\Z_{p}$-torsors of arbitrarily large essential dimension defined over finitely generated extensions of $k$.
	\end{proof}
\end{proposition}

\section{Continuous and fce dimensions}\label{ced}

The second variant of essential dimension, continuous essential dimension, is more subtle: given a pro-finite group scheme, it is defined for each $G$-torsor, while finite type essential dimension was only defined for $G$.

\subsection{Definition and first properties}

Let $G$ be a pro-finite group scheme and $T$ a $G$-torsor, consider the projective system of torsors $(T\times^GH)_{G\to H}$ where $H$ is a finite group scheme and $G\to H$ is an homomorphism. If we have two homomorphisms $G\to H$, $G\to H'$ with a third homomorphism $H\to H'$ that makes the diagram commute, a basic property of essential dimension tells us that $\ed_k T\times^GH\ge \ed_k T\times^GH'$, i.e. essential dimension increases along the projective system.

If we think the torsors $T\times^GH$ as increasingly better approximations of $T$ thanks to \autoref{proftor}, then it makes sense to consider the limit of the essential dimensions of $T\times^GH$. Thanks to the argument above, this limit exists in $\N\cup\{\infty\}$ and is just the supremum of $\ed_kT\times^GH$ where $G\to H$ varies as above. Hence, we give the following definition.

\begin{definition}\label{ceddef}
	Let $\Phi$ be a pro-finite gerbe over a field $k$, and $s:\spec L\to\Phi$ a section.
	
	The \emph{continuous essential dimension} $\ced_{k}(s)$ is the supremum of the essential dimensions of $\psi(s):\spec L\to\Phi\xar{\psi}\Psi$, where $\psi:\Phi\to\Psi$ varies among all morphisms from $\Phi$ to a finite gerbe $\Psi$.
	
	The continuous essential dimension $\ced_{k}(\Phi)$ of $\Phi$ is the supremum of $\ced_{k}(s)$, where $s$ varies among all sections $\spec K\to\Phi$ and all field extensions $K/k$. If $G$ is a pro-algebraic group scheme, we write $\ced_k G$ for $\ced_k BG$.
\end{definition}

\begin{remark}\label{whyced}
	A natural question is why we should take the limit defining the continuous essential dimension at the level of torsors and not at the level of groups, i.e. why not define $\ced_k G$ as $\lim_i\ed_kG_i$ for a pro-finite group scheme $G=\projlim_i G_i$. Obviously, this depends on taste. From our point of view, there are three reasons.
	
	\begin{itemize}
		\item It may happen that the limit $\lim_i\ed_k G_i$ does not exists, and we don't see any particular reason to prefer $\liminf_i\ed_kG_i$ or $\limsup_i\ed_kG_i$. On the other hand, the limit always exists at the level of torsors.
		\item The limit $\lim_i\ed_kG_i$ depends on the presentation of $G=\projlim_iG_i$ as a projective limit, while our definition depends only on $G$.
		\item Most importantly, we are interested in studying $G$-torsors, and $\lim_i\ed_k G_i$ depends on $G_i$-torsors that do not extend to $G$.
	\end{itemize}  
\end{remark}

Finally, we can merge in an obvious way the finite type and continuous essential dimensions and define the \emph{fce dimension} $\fced_{k} \Phi$ of a pro-finite gerbe $\Phi$.

\begin{definition}
	If $\Phi$ is a pro-finite gerbe, the \emph{fce dimension} $\fced_{k}\Phi$ of $\Phi$ is the supremum of the continuous essential dimensions $\ced_{k}s$ where $\spec L\to\Phi$ is a section over a field $L$ finitely generated over $k$. If $G$ is a pro-finite group scheme, we write $\fced_k G$ for $\fced_k BG$.
\end{definition}

\begin{lemma}\label{algebraic}
	If $\Phi$ is a finite gerbe over $k$, then $\ed_{k} \Phi=\fed_{k} \Phi=\ced_{k} \Phi=\fced_{k} \Phi$.\qed 
\end{lemma}

\begin{lemma}\label{properties}
	Let $\Phi$ be a pro-finite gerbe over a field $k$.
	\begin{enumerate}[(i)]
		\item For every section $s:\spec L\to\Phi$, we have $\ced_{k} s\le\ed_{k} s$.
		\item \[\begin{tikzcd}[row sep=tiny]
																								&	\ced_{k}\Phi\ar[dr,start anchor=-15,end anchor=165,phantom,"\le",sloped]	&					\\
					\fced_{k}\Phi\ar[ur,start anchor=15,end anchor=-165,phantom,"\le",sloped]\ar[dr,start anchor=-15,end anchor=165,phantom,"\le",sloped]		&				&	\ed_{k}\Phi		\\
																								&	\fed_{k}\Phi\ar[ur,start anchor=15,end anchor=-165,phantom,"\le",sloped]	&
				\end{tikzcd}\]
		\item If $k'/k$ is an extension, $\ced_{k'}\Phi_{k'}\le\ced_{k} \Phi$. If $k'/k$ is finitely generated, the inequality holds for $\fed$ and $\fced$ too.
	\end{enumerate}
	\begin{proof}
		The proof of (i) follows directly from the definition, (ii) follows from (i) and (iii) is identical to the analogous fact for classical essential dimension.
	\end{proof}
\end{lemma}

\begin{lemma}\label{liminfced}
	Let $\Phi=\projlim_i\Phi_i$ be a pro-finite gerbe, with $\Phi_{i}$ finite for every $i$. Then
	\[\ced_{k}\Phi\le\liminf_{i}\ed_{k}\Phi_{i}.\]
	\begin{proof}
		Since every morphism $\Phi\to\Psi$ to a finite gerbe $\Psi$ factors as $\Phi\to \Phi_{i}\to \Psi$ for some $i$, then for every section $\spec L\to\Phi$ with images $s_i:\spec L\to\Phi_i$ we have 
		\[\ced_{k}(s)=\sup_{i}\ed_{k}(s_i)=\lim_{i}\ed_{k}(s_i)\le\liminf_{i}\ed_{k}\Phi_{i}.\]
	\end{proof}
\end{lemma}

In dimension $0$, essential dimension and continuous essential dimension coincide for pro-finite gerbes.

\begin{proposition}\label{ced0}
	Let $\Phi$ be a pro-finite gerbe over $k$, and let $s:\spec L\to\Phi$ be a section where $L/k$ is an extension of fields. Then $\ced_{k}s=0$ if and only if $\ed_{k} s=0$.
	\begin{proof}
		Since $\ced_{k} s\le\ed_{k} s$, one implication is obvious. Let us assume that $\ced_{k} s=0$. Up to replacing $k$ with $\overline{k}^{L}$, we may suppose that $k$ is algebraically closed in $L$.
		
		Write $\Phi=\projlim_{i}\Phi_{i}$ with $\Phi_{i}$ finite gerbes. By hypothesis, $s_i:\spec L\to\Phi\to\Phi_i$ is defined over $k$, i.e. there exist sections $r_i:\spec k\to\Phi_i$ with $2$-commutative diagrams
		\[\begin{tikzcd}
			\spec L\rar["s"]\dar		&	\Phi\dar	\\
			\spec k\rar["r_i"]			&	\Phi_i
		\end{tikzcd}\]
		
		We want to show that the $r_{i}$ form a projective system whose limit is a section $r:\spec k\to\Phi$ such that $r_{L}\simeq s$.
		
		Let $j\ge i$ be indexes in the projective system, and define $r_{j,i}\in\Phi_i(k)$ the image of $r_j$ in $\Phi_i$. We want to give isomorphisms $r_{j,i}\simeq r_{i}$ for every $j\ge i$. Now, $\uisom_{\Phi_{i}}(r_{j,i},r_{i})$ is a finite scheme with an $L$-rational point, because we have isomorphisms
		\[r_{j,i,L}\simeq s_{i}\simeq r_{i,L}.\]
		Since $k$ is algebraically closed in $L$ and $\uisom_{\Phi_{i}}(r_{j,i},r_{i})$ is finite, the isomorphism $r_{j,i,L}\simeq r_{i,L}$ given above is defined over $k$, i.e. it is the base change of an isomorphism $\alpha_{i,j}:r_{j,i}\simeq r_{i}$.
		
		These isomorphisms respect the cocycle condition: if $j\ge i\ge h$ and we write $\phi_{h,i}:\Phi_i\to\Phi_h$, we have
		\[\alpha_{h,i}\circ\phi_{h,i}(\alpha_{i,j})=\alpha_{h,j}.\]
		In fact, this equality can be checked after base change to $L$, and over $L$ it amounts to the commutativity of the following diagram:
		\[\begin{tikzcd}[row sep=small,column sep=small]
															&	r_{i,h,L}\ar[d,"\sim" description,sloped,no head]\ar[ddr,"\alpha_{h,i}"]	&							\\[15pt]
															&	s_h																			&								\\
			r_{j,h,L}\ar[rr,"\alpha_{h,j}",swap]\ar[uur,"\phi_{h,i}(\alpha_{i,j})"]\ar[ur,"\sim" description,sloped,no head]	&			&	r_{h,L}\ar[ul,"\sim" description,sloped,no head]	
		\end{tikzcd}\]
		which is obvious. Hence $r=\projlim_{i}r_{i}:\spec k\to\Phi$ is a section, and clearly $r_{L}\simeq s$.
	\end{proof} 
\end{proposition}

\subsection{Tate modules of tori}

If $A$ is a semi-abelian variety, let us denote by $TA=\prod_{l}T_{l}A=\upi_{1}(A)$ the global Tate module, where $l$ varies among all prime numbers: it is the Nori fundamental group scheme of $A$.

\begin{lemma}\label{abband}
	If $A^1$ is a torsor for a semi-abelian variety $A$ over a field $k$, then $\Pi_{A^{1}/k}$ is banded by $TA$.
	\begin{proof}
		Kummer theory gives us an homomorphism $\H^{1}(k,A)\to\H^{2}(k,A[n])$ for every $n$, let $\Psi_{n}\in\H^{2}(k,A[n])$ be the image of $A^{1}$, it is the gerbe of liftings of $A^{1}$ along $A\xar{n}A$ and it is banded by $A[n]$. Call $S$ a copy of $A^{1}$. If $\rho:A\times A^{1}\to A$ is the action, the isomorphism 
		\[A\times S=A\times A^{1}\xar{(\rho,p_{2})} A^{1}\times A^{1}=A^{1}\times S\]
		gives us a natural trivialization of $A^{1}$ over $S$. In particular we have a preferred lifting of $A^{1}|_{S}$ along $A\xar{n}A$ and thus a preferred morphism $S\to\Psi_{n}$. These morphisms are compatible for varying $n$, thus passing to the limit we get a morphism $A^{1}=S\to\Psi=\projlim_{n}\Psi_{n}$. This in turn induces a morphism $\Pi_{A^{1}/k}\to\Psi$: it is an isomorphism, it can be checked after passing to a splitting field of $A^{1}$. Since $\Psi$ is naturally banded by $TA$, it follows that $\Pi_{A^{1}/k}$ is banded by $A$.
	\end{proof}
\end{lemma}

Thanks to \autoref{fed}, if $D$ is a torus then $\fed T_{l}D=\infty$. With continuous essential dimension we get a much more interesting result.

\begin{theorem}\label{fcedzp1}
	Let $D^{1}/k$ be a $D$-torsor over a field $k$, where $D/k$ is an algebraic torus and $l$ any prime number (possibly equal to $\on{char}k$). Then 
	\[\ced_{k}\Pi_{D^1/k,l}=\fced_k\Pi_{D^1/k,l}=\ced_{k}\Pi_{D^1/k}=\fced_k\Pi_{D^1/k}=\dim D^{1}.\]
	In particular, for $D^{1}=D$ we have
	\[\ced_{k}T_{l}D=\fced_{k}T_{l}D=\ced_{k}TD=\fced_{k}TD=\dim D.\]
	\begin{proof}
		Let $k'/k$ be a finite, separable splitting extension for $D$ (and thus $D^{1}$) and $d=\dim D$. Since $\fced_{k'}\le\fced_{k}\le\ced_{k}$, it is enough to prove $d\le\fced_{k'}\Z_{l}(1)$ and $\ced_{k}\Pi_{D^1/k}\le d$.
		
		For the lower bound, using Kummer theory it is straightforward to check that $\H^{1}(L,\Z_{l}(1)^{d})\to\H^1(L,\mu_{l}^{d})$ is surjective for any extension $L/k'$, thus $\fced\Z_{l}(1)^{d}\ge\ed\mu_{l}^d=d$.
		
		Let us prove the upper bound $\ced_{k}\Pi_{D^1/k}\le d$. Let $r=[k':k]$ be the degree of the splitting field: for any field extension $L/k$, there exists a finite separable extension $L'/L$ such that $[L':L]\le r$ and $k'/k$ is a subextension of $L'/k$ (we find $L'$ as a quotient of $L\otimes k'$). Since $\mb{G}_m^d$ has trivial cohomology, the extension map $\H^{1}(L,D)\to \H^{1}(L',D)=\{0\}$ is trivial. This implies that $\H^{1}(L,D)$ is $[L':L]$-torsion thanks to \cite[Chapter VII, \S 7, Proposition 6]{se79} (which is, "the period divides the index"). Hence, $\H^1(L,D)$ is $r!$-torsion for every extension $L/k$.
		
		Now let $L/k$ be any extension and $s:\spec L\to\Pi_{D^{1}/k}$ a section. Since $\Pi_{D'/k}$ is banded by $TD$, for every $n$ we can push the $TD$-gerbe $\Pi_{D'/k}$ along $TD\to D[n]$ and obtain a $D[n]$-gerbe $\Psi_{n}$ with an induced section $s_{n}:\spec L\to\Psi_{n}$. It is enough to prove that $\ed s_{n}\le d$ for every $n$.
		
		For every $n$, the fiber product $D^{1}_{n}=D^{1}\times_{\Psi_{n}}s_{n}$ defines a $D$-torsor such that $nD^{1}_{nm}=D^{1}_{m}\in\H^{1}(L,D)$ and $D^{1}_{1}=D^{1}_{L}$. In particular, $D^{1}_{n}\in\H^{1}(L,D)$ is divisible for every $n$. Since the group $\H^{1}(L,D)$ is $r!$-torsion, it follows that $D^{1}_{n}$ is trivial i.e. we have a section $\sigma_{n}:\spec L\to D^{1}_{n}$. By construction, the following diagram is $2$-commutative
		\[\begin{tikzcd}
			D^{1}_{n}\rar\dar	&	\spec L\dar["s_{n}"]\lar[bend right,swap,"\sigma_{n}"]	\\
			D^{1}\rar			&	\Psi_{n}
		\end{tikzcd}\]
		
		It follows that $s_{n}$ descends to the residue field of $\sigma_{n}$ composed with $D^{1}_{n}\to D^{1}$, which has transcendental degree at most $d$ over $k$.
	\end{proof}
\end{theorem}

\begin{proposition}
	For each inequality of \autoref{properties}, there is an example such that the inequality is strict.
	\begin{proof}
		For point (i), over any field $k$ with $\on{char}k\neq l$, there exists a $\Z_l(1)$-torsor $T$ with $\ed T=\infty$ thanks to \autoref{crit}, but $\ced\Z_{l}(1)=1$ and thus $\ced T\le 1$.
		
		For point (ii), if $k$ is finitely generated over $\Q$, by \autoref{crit}, \autoref{fed} and \autoref{fcedzp1} we get
		\[\fced_{k}\Z_{l}\le\fed_{k}\Z_{l}=0<1=\ced_{\bar{k}}\Z_{l}(1)=\ced_{\bar{k}}\Z_{l}\le\ced_{k}\Z_{l}<\infty=\ed_{k}\Z_{l}.\]
		
		Over any field with $\on{char}k\neq l$,
		\[\fced_{k}\Z_{l}(1)\le\ced_{k}\Z_{l}(1)=1<\infty=\fed_{k}\Z_{l}(1)\]
		and
		\[\fed_{k}\Z_{l}=0<\ed_{l}\Z_{l}.\]
		
		For point (iii), take any finite group $G$ over a field $k$ with a finite extension $k'/k$ such that $\ed_{k'}G_{k'}<\ed G$. We conclude since all variants of essential dimension coincide for $G$ finite.  
	\end{proof}
\end{proposition}

One might expect that the inequality of \autoref{liminfced} is an equality, at least when we have a presentation $G=\projlim_{i} G_{i}$ with every $G\to G_{i}$ surjective. Using \autoref{fcedzp1}, the following counterexample of F. Scavia shows that this is not the case.

\begin{example}\label{scavia}
	Let $D$ be the $1$-dimensional torus $x^2+y^2=1$ over $\Q$, it splits over $\Q(i)$. Let $T_2D=\projlim_n D[2^n]$. Thanks to \autoref{fcedzp1}, we have $\ced_{\Q}T_2D=\dim D=1$. We are now going to show that $\ed D[m]\ge 2$ for every $m$ with $4|m$. Since $D[2^{v_{2}(m)}]\s D[m]$, it is enough to show that $\ed_{\Q}D[2^n]\ge2$ for $n\ge 2$. We are going to use some concepts and results from \cite{lmmr13}.
	
	Let $\Gamma=\gal(\Q(i)/\Q)=\Z/2\Z$, the character module $M$ of $D$ is $\Z$ where $\Gamma$ acts by $x\mapsto -x$. The character module of $D[2^n]$ is $M/2^nM$. A permutation module $P$ is a $\Gamma$-module which is free as $\Z$-module, and such that $\Gamma$ acts by permutations of a basis. A $2$-presentation of $M/2^nM$ is a morphism $\phi:P\to M/2^nM$ such that $P$ is a permutation module and the cokernel is finite of odd order, in our case this is equivalent to surjectivity since $M/2^nM$ is finite of even order.
	
	Thanks to \cite[Corollary 5.1]{lmmr13}, we have that 
	\[\ed_{\Q}(D[2^n])\ge\ed_{\Q}(D[2^n];2)=\min\on{rk}\ker\phi\]
	where $\phi$ ranges among all $2$-presentations of $M/2^nM$.
	
	Let $\phi:P\to M/2^nM$ be a $2$-presentation, since $M/2^nM$ is finite we have $\on{rk}\ker\phi=\on{rk} P$. If $n>1$, the action of $\Gamma$ is non-trivial on $M/2^nM$, thus $\Gamma$ must act non-trivially on $P$ too. But then $\on{rk}P>1$, because a rank $1$ permutation module is a trivial Galois module.
\end{example}

\subsection{Fce dimension and anabelian geometry}

A. Vistoli observed that, if Grothendieck's section conjecture is true, then the étale fundamental group scheme $\upi_{1}(X)$ of a hyperbolic curve over a field finitely generated over $\Q$ should somehow have essential dimension $1$. With finite type essential dimension, we made his observation formal for proper curves. However, finite type essential dimension is still not the right tool for affine curves.

\begin{proposition}\label{nofed}
	Let $X$ be a smooth, affine curve over any field. If $\on{char}k=0$, assume $X\neq\A^1$. Then $\fed_k\Pi_{X/k}=\infty$.
	\begin{proof}
		Suppose first that $\deg\bar{X}\setminus X\ge 2$, where $\bar{X}$ is the smooth completion. Up to a finite extension of the base field we may suppose that $\bar{X}\setminus X$ has two rational points. Choose any prime $l\neq\on{char}k$. A rational point of $\bar{X}\setminus X$ induces a so-called \emph{packet of cuspidal sections} of $\Pi_{X/k}$, see \cite{sti12}. Any section $x$ of said packet induces a morphism $\Z_{l}(1)\to\upi_{1}(X,x)$, see \cite[\S 8]{bre20}. By passing to the abelianization, we check that this homomorphism is injective, thus $\fed\upi_{1}(X)=\infty$ thanks to \autoref{fed}. 
		
		If $g(X)\ge 1$, or $g(X)=0$ and $\deg\bar{X}\setminus X\ge 2$, there exists $X'\to X$ a non-trivial, connected finite étale cover. If $g(X)=0$ and $\deg\bar{X}\setminus X=1$, then $X=\A^{1}$ and by hypothesis $\on{char}k=p>0$, thus we have the Abhyankar cover $\mb{G}_m\to\A^1$ $x\mapsto x^p+1/x$, see \cite[Theorem 1]{ab57}.
		
		In any case, we have a non-trivial finite étale cover $X'\to X$ and $\deg\bar{X}'\setminus X'\ge 2$, thus $\fed_{k}\Pi_{X'/k}=\infty$. Up to a finite extension of $k$, we may suppose we have a rational point $x'\in X'(k)$ with image $x\in X(k)$, thus $\infty=\fed\upi_{1}(X',x')\le\fed\upi_{1}(X,x)$.
	\end{proof}
\end{proposition}

Because of \autoref{nofed}, the section conjecture says nothing about the finite type essential dimension of the fundamental group of an affine hyperbolic curve. This is not the case for fce dimension.

\begin{proposition}\label{dscimpl}
	Let $k$ be a field finitely generated over $\Q$, and $X$ a smooth, geometrically connected hyperbolic curve. If Grothendieck's section conjecture is true, then $\fced_{k}\Pi_{X/k}=1$.
	\begin{proof}
		If Grothendieck's section conjecture is true, then
		\[\Pi_{X/k}(k)=X(k)\sqcup\bigsqcup_{x\in\bar{X}\setminus X(k)}\mc{P}_{x}\]
		where $\mc{P}_{x}$ is, for each $x\in\bar{X}\setminus X(k)$, a so-called "packet of tangential sections", see \cite{sti12}. This is functorially isomorphic to $\H^{1}(k,\hz(1))$, see for instance \cite[\S 8]{bre20}. Since $\fced\hz(1)=1$, the statement follows.
	\end{proof}
\end{proposition}

One can read \autoref{fcedzp1} as saying that the "dimensional analogue" of the section conjecture holds for tori: if $D$ is a torus, then $\fced \Pi_{D/k}=\dim D$, just as the section conjecture predicts $\fced\Pi_{X/k}=\dim X$ if $X$ is a smooth, hyperbolic curve over a field finitely generated over $\Q$. However, for tori the dimensional statement holds over any field $k$: arithmetic plays no role, we just use Kummer theory and Hilbert's theorem 90.

For abelian varieties, we cannot expect such a result over any field: for instance, if $A$ is an abelian variety over $\C$ we have $\fced\Pi_{A/\C}=\fced\hz^{2\dim A}=2\dim A$. We are going to show that, thanks to Faltings' theorem (Tate conjecture for abelian varieties), the equality $\fced\Pi_{A/k}=\dim(A)$ holds if $k$ is a finitely generated field extension of $\Q$ and $A$ is an abelian variety over $k$.

Faltings proved that, if $A,B$ are abelian varieties over a field $k$ finitely generated over $\Q$ an $l$ is a prime number, the natural map $\hom(A,B)\otimes_{\Z}\Z_{l}\to\hom(T_{l}A,T_{l}B)$ is bijective. We need to give a "base point free" reformulation of Faltings' theorem for global Tate modules and for torsors under abelian varieties (we need this already to prove the theorem for abelian varieties, not only for torsors). 

\begin{lemma}
	Let $A,B$ be abelian varieties over a field $k$ finitely generated over $\Q$. Then
	\[\hom(TA,TB)=\widehat{\hom(A,B)}.\]
	\begin{proof}
		Since $\hom(A,B)$ is finitely generated, by Faltings' theorem we have
		\[\hom(TA,TB)=\prod_{l}\hom(T_{l}A,T_{l}B)=\prod_{l}\hom(A,B)\otimes_{\Z}\Z_{l}=\widehat{\hom(A,B)}.\]
	\end{proof}
\end{lemma}

\begin{corollary}\label{faltnice}
	Let $A,B$ be abelian varieties over a field $k$ finitely generated over $\Q$, and let $f:TA\to TB$ be an homomorphism. For every homomorphism $\phi:TB\to G$ with $G$ a finite group scheme, there exists an homomorphism $g:A\to B$ such that
	\[\phi\circ f=\phi\circ Tg.\]
	The same statement holds for $T_{l}\bullet$ for any prime $l$.\qed
\end{corollary}

\begin{lemma}[Base point free Faltings' theorem]\label{faltbpf}
	Let $k$ be a field finitely generated over $\Q$, and $E,F\to\spec k$ torsors for abelian varieties $A,B$ over $k$ and $l$ a prime number. Let $\Pi_{E/k},\Pi_{F/k}$ be the étale fundamental gerbes of $E,F$, and $\rho:\Pi_{E/k}\to\Pi_{F/k}$ a morphism. 
	
	For every finite gerbe $\Phi$ and every morphism $\phi:\Pi_{F/k}\to\Phi$ there exists a $B$-torsor $F'$ and morphisms $f:E\to F'$, $\Pi_{F'/k}\to\phi'$ such that the following diagram $2$-commutes
	\[\begin{tikzcd}
		E\rar["\pi_E"]\dar["f"]		&	\Pi_{E/k}\rar["\rho"]\dar["\pi(f)"]		&	\Pi_{F/k}\dar["\phi"]		\\
		F'\rar["\pi_{F'}"]			&	\Pi_{F'/k}\rar["\phi'"]					&	\Phi
	\end{tikzcd}\]
	
	The same statement holds for $\Pi_{\bullet/k,l}$ for any prime $l$.
	\begin{proof}
		Let us prove the statement for $\Pi_{\_/k}$, the one for $\Pi_{\_/k,l}$ is analogous. Thanks to \cite[Lemma 5.12]{bv15}, we may suppose that $\Pi_{F/k}\to\Phi$ is Nori-reduced (this essentially amounts to the fact that the induced homomorphism between automorphism groups is surjective), and thus we may suppose that $\Phi$ is abelian since $\Pi_{F/k}$ is abelian. 
		
		Thanks to \autoref{abband}, the bands of $\Pi_{E/k}$ and $\Pi_{F/k}$ are respectively $TA$ and $TB$, call $G$ the band of $\Phi$. We have induced morphisms of bands $TA\xar{\rho} TB\xar{\phi} G$. Since $G$ is finite, by \autoref{faltnice} we have a morphism of abelian varieties $f:A\to B$ such that 
		\[\phi\circ Tf=\phi\circ\rho.\]
		
		Write 
		\[F'=E\times^AB\]
		the induced torsor along $f:A\to B$, we have a natural equivariant morphism $E\to F'$. By construction, we have a natural isomorphism
		\[\Pi_{F'/k}\simeq\Pi_{E/k}\times^{Tf}TB\]
		between the fundamental gerbe of $F'$ and the induced gerbe along $Tf$. Since $\phi\circ TF=\phi\circ\rho$ the diagram
		\[\begin{tikzcd}
			TA\rar["\rho"]\dar["Tf"]	&	T B\dar["\phi"]		\\
			TB\rar["\phi"]				&	G
		\end{tikzcd}\]
		commutes. If we push $\Pi_{E/k}$ along this diagram we get the commutative diagram
		\[\begin{tikzcd}
			\Pi_{E/k}\rar["\rho"]\dar["\pi(f)"]		&	\Pi_{F/k}\dar["\phi"]		\\
			\Pi_{F'/k}\rar["\phi'"]					&	\Phi
		\end{tikzcd}\]
		as desired.
	\end{proof}
\end{lemma}

\begin{lemma}
	Let $V$ be a smooth variety over a field $k$ of characteristic $0$ with with Albanese torsor $V\to A^1$, and let $\Pi^{\rm ab}_{V/k}$ be the abelianized fundamental gerbe of $V$. The natural morphism
	\[\Pi^{\rm ab}_{V/k}\to \Pi_{A^1/k}\]
	is a quotient of gerbes, and the kernel is torsion.
	\begin{proof}
		This is classical for $k$ algebraically closed and étale fundamental groups, see \cite[Corollary 5.8.10]{sz09}. The general case follows from the fact that the étale fundamental gerbe behaves well under base change, see \cite[Proposition 6.1]{bv15}.
	\end{proof}
\end{lemma}

\begin{corollary}\label{abelianized}
	If $V$ is a smooth projective variety over a field $k$ of characteristic $0$ with Albanese torsor $V\to A^1$ and $V\to\Phi$ is a morphism to an abelian, torsion-free gerbe $\Phi$ then we have a factorization $V\to A^1\to\Phi$.\qed
\end{corollary}

\begin{lemma}\label{lowab}
	Let $A$ be an abelian variety over a field $k$ and $n>1$ an integer. Let $T\to\spec k(A)$ the restriction to the generic point of the $A[n]$-torsor $A\to A$ induced by multiplication by $n$. Then $\ed T=\dim A$.
	\begin{proof}
		Clearly $\ed T\le\dim A$, we want to prove $\dim A\le\ed T$. Kummer theory gives us an exact sequence
		\[A(k(A))\xar{n}A(k(A))\to\H^{1}(k(A),A[n])\to\H^{1}(k(A),A).\]
		
		Let $k(A)/L/k$ be a subextension such that $T$ descends to a $A[n]$-torsor $Q\to L$. Let $A^{1}\in\H^{1}(L,A)$ be the image of $Q\in\H^{1}(L,A[n])$, it is an $A$-torsor over $L$. Since $Q_{k(A)}=T$ comes from $A(k(A))$, it follows that $k(A)$ splits $A^{1}$. In particular, we have a point $\spec L(A_{L})\to A^{1}$. Since $A_{L}$ has a smooth $L$-rational point, it follows that $A^{1}$ has an $L$-rational point too by the Lang-Nishimura theorem. It is possible to avoid the Lang-Nishimura theorem: morphisms to abelian varieties always extend, and by Galois descent this generalizes to abelian torsors, thus we have a morphism $A_{L}\to A^{1}$.
		
		It follows that $Q$ comes from an $L$-rational point $\spec L\to A$. The point $\spec L\to A$ extends to a morphism $s:V\to A$ where $V$ is some regular variety with $L=k(V)$, while the field extension $k(A)/L$ gives a dominant rational morphism $f:A\dashrightarrow V$. The composition $s\circ f$ is a morphism $A\to A$ since morphisms to abelian varieties always extend. We have that $s\circ f(A)=a+B\s A$ for some $a\in A(k)$ and $B\s A$ a sub-abelian variety. The Kummer exact sequence applied to the isomorphism $Q_{k(A)}\simeq T$ gives us
		\[\id-s\circ f\in n\on{Mor}(A,A)\s\on{Mor}(A,A)\]
		since $Q_{k(A)}$ is associated to $\spec k(A)\to A\xar{s\circ f}A$ while $T$ is associated to the generic point $\spec k(A)\to A$. The morphism $\id-s\circ f: A\to A$ induces the translation $\tau_{-[a]}$ by $-[a]$ on $A/B$, thus
		\[\tau_{-[a]}\in n\on{Mor}(A/B, A/B).\]
		Since a translation is an isomorphism and $n>1$, it follows that $A/B$ is trivial and $A=B$, i.e. $s\circ f(A)=A$. In particular, $s$ is dominant, i.e. $\spec L\to A$ maps to the generic point and hence $\trdeg(L/k)\ge \dim A$.
	\end{proof}
\end{lemma}

\begin{theorem}\label{dscab}
	Let $A^1$ be an $A$-torsor for an abelian variety $A$ over a field $k$ finitely generated over $\Q$, and $l$ a prime number. Then 
	\[\fced_k\Pi_{A^1/k}=\fced_k\Pi_{A^1/k,l}=\dim A^1.\]
	In particular, for $A^1=A$ we have
	\[\fced_kTA=\fced_kT_lA=\dim A.\]
	\begin{proof}
		Let us first prove the lower bound. Up to a finite extension of $k$, we may suppose that $A^1=A$ is an abelian variety, and we focus on a prime $l$. By \autoref{lowab} the $A[l]$-torsor on $\spec k(A)$ given by multiplication by $l$ has essential dimension $\dim A$. Since this lifts to a $T_{l}A$-torsor on $\spec k(A)$, we get the lower bound.
	
		Let us now prove the upper bound. We do it for $\Pi_{A^1/k}$, the argument for $\Pi_{A^1/k,l}$ is analogous.
		
		Let $k'/k$ be a field finitely generated over $k$, and $\spec k'\to\Pi_{A^1/k}$ a section. Up to replacing $k$ with $\bar{k}^{k'}$, we may suppose that $k$ is algebraically closed in $k'$. By resolution of singularities there exists a smooth, geometrically connected projective variety $V$ with $k(V)=k'$. For every prime $l$, since $\Z_{l}(1)$ has weight $-2$ and $T_{l}A$ has weight $-1$, there are no non-trivial homomorphisms $\Z_{l}(1)\to T_{l}A$, and this remains true after finitely generated extensions. Hence, thanks to \autoref{r-proper}, $\spec k'\to\Pi_{A^1/k}$ extends uniquely to a morphism $V\to\Pi_{A^1/k}$.
		
		Let $V\to B^1$ be the Albanese torsor of $V$, it is a torsor for the Albanese variety $B$. Since $\Pi_{A^1/k}$ is abelian and torsion free, by \autoref{abelianized} we have a factorization
		\[V\to B^1\to\Pi_{B^1/k}\to \Pi_{A^1/k}.\]
		
		Let us suppose we have a morphism $\phi:\Pi_{A^1/k}\to\Phi$ with $\Phi$ a finite gerbe, we have to show that the composition $\spec k'\to V\to\Pi_{A^1/k}\to\Phi$ factorizes through a field of transcendence degree less than or equal to $\dim A$.
		
		By \autoref{faltbpf} there exists a morphism $f:B\to A$ such that, if $A^{f}=B^1\times^BA$ is the induced $A$-torsor, the following diagram commutes:
		\[\begin{tikzcd}
			V\rar	&	B^1\rar\dar["f"]	&	\Pi_{B^1/k}\rar\dar["\pi(f)"]	&	\Pi_{A^1/k}\dar["\phi"]		\\
					&	A^{f}\rar			&	\Pi_{A^{f}/k}\rar["\phi"]		&	\Phi
		\end{tikzcd}\]
		In particular, this tells us that the composed morphism $\spec k'\to\Phi$ factorizes through the residue field of a point of $A^{f}$, which has transcendence degree less than or equal to $\dim A$, as desired.
	\end{proof}
\end{theorem}

\appendix

\section{Extension of torsors}\label{extapp}

Finite type essential dimension focuses on torsors defined over fields finitely generated over the base field. In this context, an often meaningful question is the following: if we have a torsor defined on the generic point of a variety, does it extend to the whole variety? In this appendix, which is completely independent from the concept of essential dimension, we prove some results regarding this question.

We are going to use extensively the notion of gerbes and in particular we are going to replace étale fundamental groups with étale fundamental gerbes, see \cite{bv15}. This is not strictly necessary, but it allows us to handle better various situations where fixing a base point is troublesome. Every result about gerbes can be translated into a result about torsors by considering the gerbe $BG$ if $G$ is a group scheme. 

\subsection{Maximal open subset of definition}

\begin{lemma}\label{extension}
	Let $G$ be a pro-finite étale group scheme over a field $k$, $V$ a geometrically integral scheme and $T=\spec A\to\spec k(V)$ a torsor defined on the generic point. Let $\tilde{A}\s A$ be the normalization of $\O_V$ in $A$, and $\tilde{T}\to V$ the relative spectrum.
	
	\begin{itemize}
		\item The action of $G$ on $T$ extends to $\tilde{T}$.
		\item If $\tilde{T}\to V$ is étale, then $V$ is normal and $\tilde{T}\to V$ is a $G$-torsor. 
		\item On the other hand, if $V$ is normal and an extension of $T$ to $V$ exists, then it is unique and it coincides with $\tilde{T}$.
	\end{itemize}
	\begin{proof}
		The problem is local, we may suppose $V=\spec\tilde{B}$ where $\tilde{B}$ is an integral $k$-algebra with fraction field $k(V)=B$. Moreover, it is straightforward to get the general case from the one in which $G$ is finite étale, hence we make this assumption.
		
		By definition, $\tilde{A}$ is the integral closure of $\tilde{B}\s B$ in $A$. Using the Yoneda lemma, $G$ acts by ring homomorphisms on $A$ and fixes the elements of $\tilde{B}\s B\s A$, hence elements integral over $\tilde{B}$ are sent to elements integral over $\tilde{B}$ by the action, i.e. the action of $G$ restricts to $\tilde{A}$. Here we are subtly using the fact that $G$ is étale and integral closure commutes with étale base change, see \cite[Proposition 18.12.15]{egaiv4}. For a more down-to-earth proof, write $G=\spec R$ and consider the co-module structure
		\[A\to A\otimes R.\]
		The elements of $\tilde{A}$ are integral over $\tilde{B}$ and thus they are mapped to elements of $A\otimes R$ integral over $\tilde{B}\otimes R$. Since integral closure commutes with étale base change, the elements of $A\otimes R$ integral over $\tilde{B}\otimes R$ are precisely the elements of $\tilde{A}\otimes R$.
		
		Now suppose that $\tilde{T}\to V$ is étale. In particular, if $\tilde{V}$ is the normalization of $V$ in $k(V)$, we have a factorization $\tilde{T}\to \tilde{V}\to V$, hence the normalization $\tilde{V}\to V$ is étale and thus an isomorphism.
		
		Now, let $\rho:G\times \tilde{T}\to\tilde{T}$ the action. We have a natural morphism $\rho\times p_2:G\times \tilde{T}\to\tilde{T}\times_{V}\tilde{T}$: the fact that $\tilde{T}$ is a $G$-torsor is equivalent to the fact that $\rho\times p_2$ is an isomorphism. But $G\times\tilde{T}$ and $\tilde{T}\times_{V}\tilde{T}$ are both finite étale covers of $V$ and $\rho\times p_2$ is a morphism of covers, thus in order to prove that it is an isomorphism it is enough to show that it is generically an isomorphism, and this is true by hypothesis since $T$ is a torsor.
		
		On the other hand, suppose that $\tilde{T}'\to V$ is some extension of $T$ and that $V$ is normal. Then $\tilde{T}'\to V$ is finite étale (because $G$ is finite étale), hence $\tilde{T}'$ is normal and finite over $V$: the former says that $\O_{T'}\s A$ contains the normalization $\tilde{A}$ of $\tilde{B}$ (since $\O_{V}\supseteq \tilde{B}$ and is normal), the latter that it is contained in it.
	\end{proof}
\end{lemma}

\begin{corollary}\label{gextension}
	Let $\Phi$ be a finite étale gerbe over a field $k$, $V$ a normal scheme over $k$ and $s:\spec k(V)\to\Phi$ a morphism. Call $\xi$ the generic point. Then there exists an open subset $U_{\rm max}\s V$ with a morphism $u_{\rm max}:U_{\rm max}\to\Phi$ and an isomorphism
	\[\phi_{\rm max}:s\to u_{{\rm max},\xi}\]
	in $\Phi(k(V))$ such that
	\begin{enumerate}[(i)]
		\item for every other $u:U\to\Phi$, $\phi:s\to u_{\xi}$ as above, we have $U\s U_{\rm max}$ and there exists a unique isomorphism $\psi:u_{\rm max}|_{U}\to u$ such that $\psi_{\xi}\circ \phi_{\rm max}=\phi$,
		\item a point $v\in V$ is in $U_{\rm max}$ if and only if $s$ extends to $\spec \O_{V,v}$,
		\item if $V$ is regular then $V\setminus U_{\rm max}\s V$ has pure codimension $1$,
		\item if $k'/k$ is a separable extension, then $U_{\rm max}\times k'=(U\times k')_{\rm max}$ and $u_{\rm max}\times k'=(u\times k')_{\rm max}$.
	\end{enumerate}
	\begin{proof}
		If $\Phi$ has a section $\spec k\to\Phi$, points (i), (ii) and (iii) are a direct consequence of \autoref{extension}: there exists the greatest open subset $U_{\rm max}$ where a torsor $T\to\spec k(V)$ extends, the extension is unique, if $V$ is regular then $V\setminus U_{\rm max}$ has pure codimension $1$ by purity of branch locus, a point is in $U_{\rm max}$ if and only if $T$ is unramified over $v$ if and only if it extends to $\O_{V,v}$. Since integral closure commutes with étale base change (see \cite[Proposition 18.12.15]{egaiv4}), we get (iv).
		
		Otherwise, since $\Phi$ is finite étale there exists a finite Galois extension $k'/k$ with a section $\spec k'\to\Phi$, this lets us identify $\Phi_{k'}=BG'$ for some finite étale group scheme $G'$ over $k'$. Write $V'=V\times\spec k'$, $V'$ is still normal because integral closure commutes with étale base change. We have a $G'$-torsor $T'\to\spec k(V')$. For every element $\sigma\in\gal(k'/k)$ we have a $\sigma$-equivariant isomorphism of $k$-schemes $\sigma^* T'\to T'$, thus the étale loci of $T'$ and $\sigma^* T'$ coincide i.e. the étale locus $U_{\rm max}'\s V'$ of $T'$ is Galois-invariant, let $U_{\rm max}\s V$ be its image in $V$. We have that $\tilde{T'}|_{U_{\rm max}'}$ defines a morphism $U_{\rm max}'\to\Phi_{k'}=BG'$ which is Galois invariant. To check that it descends to a morphism $U_{\rm max}\to\Phi$ we only have to check the cocycle condition: but this can be checked on the generic point, where it is obviously satisfied.
		
		It is immediate to check that $U_{\rm max}\to\Phi$ satisfies the requested conditions since $U_{\rm max}'\to\Phi_{k'}$ satisfies them.
	\end{proof}
\end{corollary}

\begin{corollary}\label{gprextension}
	Let $\Phi=\projlim_i\Phi_i$ be a pro-finite étale gerbe over a field $k$, $V$ a normal scheme and $s:\spec k(V)\to\Phi$ a morphism.
	
	Write $U_{{\rm max}, i}\s V$ for the open subset given by \autoref{gextension} with respect to $\spec k(V)\to\Phi_i$, we have $U_{{\rm max}, j}\s U_{{\rm max}, i}$ if $j\ge i$. Let $U_{{\rm max}}=\bigcap_{i}U_{{\rm max},i}$.
	
	Suppose that $U_{{\rm max}}=U_{{\rm max},i}$ for $i>>0$. Then the conclusions of \autoref{gextension} hold for $\spec k(V)\to\Phi$. \qed 
\end{corollary}

\subsection{Extension to points of codimension $1$}

Let $V$ be a regular, geometrically connected scheme over $k$, $D\in V$ an irreducible Cartier divisor. Suppose we have a pro-finite étale gerbe $\Phi$ and a morphism $V\setminus D\to\Phi$, we want to study conditions under which it extends to a morphism $V\to\Phi$. 

\begin{lemma}\label{dvrext}
	Let $V\setminus D\to\Phi$ be as above. Let $\xi$ be the generic point of $D$, and suppose that $\spec k(V)\to\Phi$ extends to $\spec\O_{V,\xi}$. Then $V\setminus D\to\Phi$ extends to $V$.
	\begin{proof}
		Write $\Phi=\projlim_i\Phi_i$, and let $U_{{\rm max},i}\s V$ be as in \autoref{gprextension}. By hypothesis, $U_{{\rm max},i}$ contains $V\setminus D$. Since $V\setminus U_{{\rm max},i}$ has pure codimension $1$, we have that $U_{{\rm max},i}$ is equal either to $V\setminus D$ or $V$. By \autoref{gextension}.(ii), $\xi\in U_{{\rm max},i}$, thus $U_{{\rm max},i}=V$. The claim follows.
	\end{proof}
\end{lemma}

Thanks to \autoref{dvrext}, we have reduced the problem to the case in which $V=\spec R$ is the spectrum of a geometrically regular, geometrically connected DVR $R$ and $D=p\in \spec R$ is the closed point. This allows us to avoid problems arising from the fact that $D$ might be singular.

Let $\sqrt[\infty]{(\spec R,p)}$ be the infinite root stack of $\spec R$ at $p$, see \cite{tv18}. There are natural morphisms
\[\spec k(R)\hookrightarrow\sqrt[\infty]{(\spec R,p)}\twoheadrightarrow \spec R.\]

\begin{lemma}\label{tameroot}
	The natural morphism $\Pi^{\et}_{k(R)/k}\to\Pi^{\et}_{\sqrt[\infty]{(\spec R,p)}/k}$ is a quotient of gerbes.
	
	A morphism $\Pi^{\et}_{k(R)/k}\to\Phi$ where $\Phi$ is a finite étale gerbe factorizes through $\Pi^{\et}_{\sqrt[\infty]{(\spec R,p)}/k}$ if an only if, for any section $\spec k'\to\Phi$, the induced covering $Y\to\spec k'(R_{k'})$ is tamely ramified on closed points of $\spec R_{k'}$.
	
	In particular, $\Pi^{\et}_{k(R)/k}\simeq\Pi^{\et}_{\sqrt[\infty]{(\spec R,p)}/k}$ if $\on{char}k=0$.
	\begin{proof}
		This follows from \cite[Proposition 3.2.2]{bor09}.
	\end{proof}
\end{lemma}

If $G$ is a pro-finite group, recall that the order of $G$ is defined as a supernatural number, i.e. a formal product $\prod_{l}p^{n_{l}}$ where $l$ ranges among all prime numbers and $n_{l}\in\N\cup\{\infty\}$.

\begin{definition}
	Let $\Phi$ be a pro-finite étale gerbe over $k$. The \emph{order} $|\Phi|$ of $\Phi$ is the order of $\aut_{\Phi}(s)$ where $s\in\Phi(\bar{k})$ is any geometric section.
\end{definition}

\begin{remark}
	If $\on{char}k$ does not divide the order of a finite étale gerbe $\Phi$ and $\spec k(R)\to\Phi$ is a section, the induced morphism $\Pi^{\et}_{k(R)/k}\to\Phi$ automatically satisfies the condition of \autoref{tameroot} since ramification indices divide the order of $\Phi$. Hence, if $\Phi$ is finite étale and $\on{char}k$ does not divide $|\Phi|$, a section $\spec k(R)\to\Phi$ induces a morphism $\Pi^{\et}_{\sqrt[\infty]{(\spec R,p)}/k}\to\Phi$. Passing to the limit, this fact holds for $\Phi$ pro-finite étale, too.
\end{remark}

The fiber of 
\[\sqrt[\infty]{(\spec R,p)}\to\spec R\]
over the closed point $p$ is non-canonically isomorphic to $B_{k(p)}\hz(1)$.

\begin{definition}\label{holedef}
	We call the composition
	\[B_{k(p)}\hz(1)\to\sqrt[\infty]{(\spec R,p)}\to\Pi^{\et}_{\sqrt[\infty]{(\spec R,p)}/k}\]
	the \emph{hole} at $p$.
	
	More generally, if $V,D$ are as above and $\xi$ is the generic point of $D$, we call the composition
	\[B_{k(D)}\hz(1)\to\sqrt[\infty]{(\spec \O_{V,D},\xi)}\to\Pi^{\et}_{\sqrt[\infty]{(\spec \O_{V,D},\xi)}/k}\to\Pi^{\et}_{\sqrt[\infty]{(V,D)}/k}\]
	the \emph{hole} at $D$.
\end{definition}

Recall that a morphism $X\to Y$ of fibered categories over a field $k$ is constant if there exists a factorization $X\to\spec k\to Y$.

\begin{lemma}\label{divext}
	Let $V$ be a geometrically regular, geometrically connected scheme over a field $k$ and $D\s V$ a codimension $1$ sub-scheme.
	
	Suppose we have a pro-finite étale gerbe $\Phi$ with $\on{char}k\nmid|\Phi|$ and a morphism $V\setminus D\to\Phi$. Then $V\setminus D\to\Phi$ extends uniquely to a morphism $V\to\Phi$ if and only if the morphism $B_{k(D)}\hz(1)\to\Phi_{k(D)}$ induced by the hole at $D$ is constant.
	\begin{proof}
		The "only if" part is obvious, let us prove the other implication. Thanks to \autoref{dvrext} we may reduce to the case in which $V$ is the spectrum of a DVR $R$ and $D=p\in\spec R$ is the closed point.
		
		Write $\Phi=\projlim\Phi_i$, we may apply \autoref{gextension} and the claim is equivalent to the fact that $p\in U_{{\rm max},i}$ for every $i$.
		
		Fix an index $i$. There exists a finite separable extension $k'/k$ with a section $\spec k'\to\Phi_i$: the point $p$ may split in the extension $k'$, but it is immediate to check that the hypothesis is still satisfied for every point of $\spec R\otimes_k k'$ over $p$. Hence we may base change to $k'$ since $U_{\rm max}$ behaves well with respect to base change, and we may assume that $\Phi_i=BG$ for some finite étale group scheme $G$ with $\on{char}k\nmid|G|$.
		
		Let $T\to\spec R$ be the ramified $G$-covering associated to $\spec k(R)\to BG$. Then $[T/G]\simeq\sqrt[r]{(\spec R,p)}$ where $r$ is the ramification index of $T$ over $p$, see \cite[Lemme 3.3.1]{bor09}. By construction, 
		\[B_{k(p)}\mu_r\to\sqrt[r]{(\spec R,p)}\to BG\] is faithful. By hypothesis, $B_{k(p)}\hz(1)\to BG$ is constant. This implies that $B_{k(p)}\hz(1)\to B_{k(p)}\mu_r$ is constant too, i.e. $r=1$. This means precisely that $T$ is unramified at $p$, i.e. $p\in U_{{\rm max},i}$.
	\end{proof}
\end{lemma}

Whether a morphism is constant can be checked after a separable base change. Recall that a fibered category $X$ is concentrated if there exists an quasi-compact scheme $U$ and a representable, quasi-compact, quasi-separated, faithfully flat morphism $U\to X$. It can be easily checked that algebraic stacks of finite type over a field and pro-finite gerbes are concentrated. See \cite[Appendix A]{bre20} for the definition of geometrically connected fibered categories.

\begin{lemma}\label{constbc}
	Let $X$ be an inflexible (resp. geometrically connected) fibered category over $k$, and let $\Phi$ be a pro-finite (resp. pro-finite étale) gerbe. Let $k'/k$ be an algebraic and separable extension. Suppose that either $k'/k$ is finite, or $X$ is concentrated. A morphism $X\to\Phi$ is constant if and only if $X_{k'}\to\Phi_{k'}$ is constant.
	\begin{proof}
		The "only if" part is obvious. Suppose that $X_{k'}\to\Phi_{k'}$ is constant. We do the case in which $X$ is inflexible and $\Phi$ is pro-finite, $X$ geometrically connected and $\Phi$ pro-finite étale is analogous if we replace the results of \cite{bv15} with those of \cite[Appendix A]{bre20}.
		
		Let us first do the case in which $\Phi$ is finite. Let $X\to\Delta\to\Phi$ be such that $X\to\Delta$ is Nori-reduced and $\Delta\to\Phi$ is faithful, this factorization exists thanks to \cite[Lemma 5.12]{bv15}. We have a further factorization $X\to\Pi_{X/k}\to\Delta\to\Phi$ through the Nori fundamental gerbe, and $\Pi_{X/k}\to\Delta$ is Nori-reduced. Thanks to \cite[Proposition 6.1]{bv15}, $\Pi_{X_{k'}/k'}=\Pi_{X/k}\times_{k}k'$. Thanks to the characterization given in \cite[Proposition 3.10]{bv19}, Nori-reduced morphisms of gerbes are stable under base change, and thus $X_{k'}\to\Pi_{X_{k'}/k'}\to\Delta_{k'}$ is Nori-reduced. Since $X_{k'}\to\Phi_{k'}$ is constant, by the uniqueness part of \cite[Lemma 5.12]{bv15} we get $\Delta_{k'}=\spec k'$, thus $\Delta=\spec k$.
		
		If $\Phi$ is pro-finite, write $\Phi=\projlim_{i}\Phi_{i}$ with $\Phi_{i}$ finite. Each $X\to\Phi_{i}$ is constant by the preceding case, and the factorization through $\spec k$ is unique thanks to \cite[Lemma 5.12]{bv15}, thus we may pass to the limit. 
	\end{proof}
\end{lemma}

\begin{corollary}\label{r-proper}
	Let $V$ be a connected, geometrically regular scheme over $k$. Assume moreover that $k(p)/k$ is finitely generated for every $p\in V$.
	
	Let $\Phi$ a pro-finite étale gerbe over $k$ with $\on{char}k\nmid\Phi$, and suppose that for every finitely generated extension $k'/k$ and for every section $s:\spec k'\to\Phi$ there are no non-trivial homomorphisms $\hz(1)\to \uaut_{\Phi}(s)$ of group schemes over $k'$.
	
	Every generic section $\spec k(V)\to\Phi$ extends uniquely to a morphism $V\to\Phi$.
	\begin{proof}
		Fix a section $s:\spec k(V)\to\Phi$. Since $k(V)$ is finitely generated over $k$, the separable closure of $k$ in $k(V)$ is finite, let $k'/k$ be a Galois closure. The base change $V_{k'}$ is a disjoint union of geometrically regular, geometrically connected schemes over $k'$. Since the locus of definition of commutes with base change along finite, separable extensions of $k$ (see \autoref{gextension}), we may assume $k'=k$, i.e. $V$ geometrically connected.
		
		Write $\Phi=\projlim_i\Phi_i$, it is enough to show that every codimension $1$ point of $V$ belongs to $U_{{\rm max},i}$ for every $i$, i.e. that $\spec k(V)\to\Phi$ extends to $\spec\O_{V,\xi}$ for every codimension $1$ point $\xi$.
		
		Hence, we may fix $\xi$ and suppose that $V=\spec\O_{V,\xi}$ is the spectrum of a DVR. Then the claim follows from \autoref{divext} since the residue field of $R$ is finitely generated over $k$, and thus we may apply the hypothesis.
	\end{proof}
\end{corollary}

\subsection{Extension by dominance}

\begin{lemma}\label{smoothext}
	Let $W\to V$ be a surjective smooth morphism of normal schemes, $\Phi$ a pro-finite étale gerbe, $s_V:\spec k(V)\to\Phi$ a section and $s_W$ its composition with $\spec k(W)\to\spec k(V)$.
	
	Suppose that $s_W$ extends to a morphism $W\to\Phi$. Then $s_V$ extends to a morphism $V\to\Phi$.
	\begin{proof}
		Let $k'/k$ a separable extension with a section $\spec k'\to\Phi$. Over $k'$, the claim follows directly from \autoref{extension} plus the fact that integral closure commutes with smooth base change, see \cite[\href{https://stacks.math.columbia.edu/tag/03GG}{Lemma 03GG}]{stacks-project}. Since the maximal locus of definition commutes with base change along separable extensions (see \autoref{gprextension}) and $\spec k(V\times k')\to\Phi_{k'}$ extends to $V\times k'\to\Phi_{k'}$ we get the claim.
	\end{proof}
\end{lemma}

If we assume that $\Phi$ is torsion free and $\on{char}k\nmid|\Phi|$, we may drop the smoothness assumption in \autoref{smoothext}.

\begin{lemma}\label{ramext}
	Let $W\to V$ be a surjective morphism of connected, geometrically regular schemes over a field $k$, $\Phi$ a torsion-free pro-finite étale gerbe with $\on{char}k\nmid|\Phi|$, $s_V:\spec k(V)\to\Phi$ a section and $s_W$ its composition with $\spec k(W)\to\spec k(V)$.
	
	Suppose that $s_W$ extends to a morphism $W\to\Phi$. Then $s_V$ extends to a morphism $V\to\Phi$.
	\begin{proof}
		It is enough to show that $s_{v}$ is defined on every codimension $1$ point, i.e. we may assume that $V=\spec R$ is the spectrum of a geometrically regular DVR with closed point $p\in V$. We may also replace $W$ with the spectrum $\spec S\s W$ of a geometrically regular DVR $S$ such that $\spec S\s W\to V=\spec R$ is surjective. Let $q\in\spec S$ be the closed point.
		
		We want to apply \autoref{divext}. We want to show that the hole
		\[B_{k(p)}\hz(1)\to\Pi^{\et}_{\sqrt[\infty]{V,p}}\to \Phi\]
		factorizes through some section $\spec k(p)\to \Phi$. Call $\phi$ the image of the tautological section $\spec k(p)\to B_{k(p)}\hz(1)$ in $\Phi$. We want to show that the induced homomorphism
		\[\sigma_{V,p}:\hz(1)\to \uaut_{\Phi_{k(p)}}(\phi)\]
		is trivial. By hypothesis, this is true over $W$, i.e. the analogous homomorphism
		\[\sigma_{W,q}:\hz(1)\to \uaut_{\Phi_{k(q)}}(\phi)\]
		is trivial. But now we have a commutative diagram
		\[\begin{tikzcd}
			\hz(1)\ar[rr,"r"]\ar[dr,"\sigma_{W,c}=0",swap]	&		&	\hz(1)\ar[dl,"\sigma_{V,d}"]	\\
											&	\uaut_{\Phi_{k(c)}}(\phi)	&
		\end{tikzcd}\]
		where the horizontal arrow $\hz(1)\to \hz(1)$ is just multiplication by the ramification index $r$ of $W\to V$ at $q$. Hence, since $\sigma_{W,c}$ is trivial and $\uaut_{\Phi_{k(d)}}(\phi)$ is torsion free by hypothesis, we get that $\sigma_{V,d}$ is trivial too, as desired.
	\end{proof}
\end{lemma}

\begin{lemma}\label{ratext}
	Let $W\dashrightarrow V$ a rational dominant map of smooth projective varieties over a field $k$, $\Phi$ a torsion-free pro-finite étale gerbe with $\on{char}k\nmid |\Phi|$, $s_V:\spec k(V)\to\Phi$ a section and $s_W$ its composition with $\spec k(W)\to\spec k(V)$.
	
	Suppose that $s_W$ extends to a morphism $W\to\Phi$. Then $s_V$ extends to a morphism $V\to\Phi$.
	\begin{proof}
		There exists a smooth projective variety $W'$ with surjective morphisms $W'\to W$, $W'\to V$ which commute with the given rational map $W\dashrightarrow V$: up to replacing $W$ with $W'$, we may suppose that the rational map $W\dashrightarrow V$ is a projective dominant morphism, hence surjective. Now apply \autoref{ramext}.
	\end{proof}
\end{lemma}

\printbibliography
	
\end{document}